# Standard posets and integral weight bases for symmetric powers of minuscule representations[*]

Michael C. Strayer
Hampden–Sydney College
Hampden–Sydney, VA 23943 U.S.A.
mstrayer@hsc.edu

June 16, 2025

## Abstract

This paper extends our earlier work where we constructed "minuscule" representations of Kac–Moody algebras from colored posets in a way that maintains key properties of the well-known minuscule representations of simple Lie algebras. In this paper we work only with finite posets. We define standard posets here as ones that can be used to construct weight bases of $m^{\text{th}}$ symmetric powers ($m \geq 1$) of these minuscule Kac–Moody representations over the integers in a certain fashion. Our main result is to show that our "$\Gamma$-colored $d$-complete" and "$\Gamma$-colored minuscule" posets are standard. When the algebra at hand is a simply laced simple Lie algebra and the representation minuscule in the classic sense (i.e. isomorphic to irreducible $V(\lambda)$ for minuscule highest weight $\lambda$), our result produces a concrete combinatorially described weight basis for the irreducible representation $V(m\lambda)$ that is indexed in a natural fashion by $m$-multichains in the weight lattice for $V(\lambda)$. C.S. Seshadri first showed such an indexing of a basis is possible. Our work here is entirely combinatorial and does not use results or techniques from algebraic geometry. Constructions in this paper are independent of Lie type and actions of Kac–Moody algebra elements on basis vectors are effectively specified.

## 1  Introduction

This paper presents results jointly obtained with Robert A. Proctor. A minuscule representation of a simple Lie algebra is an irreducible highest weight representation of that algebra whose weights are all in the Weyl group orbit of the highest weight. These representations are indexed by the list of their associated highest weights, the minuscule fundamental weights [Bou, Ch. VIII, §7.3]. We expanded this notion in [Str1, Str2, Str3] to various kinds of "minuscule" representations of Kac–Moody algebras and their positive Borel subalgebras that were explicitly constructed using partially ordered sets that were colored by the nodes of the Dynkin diagrams associated to the algebras at hand. These new notions of minuscule include highest weight and non-highest weight representations, finite dimensional and infinite dimensional representations, and representations of Kac–Moody algebras and representations of their positive Borel subalgebras. These new representations retain the key way that algebra generators act on weight vectors in the classic minuscule representations of simple Lie algebras. The main result of [Str1] characterized which combinatorial properties of these "Γ-colored" posets were necessary and sufficient for these representations to be produced, resulting in new axiomatic definitions of "Γ-colored $d$-complete" and "Γ-colored minuscule" posets. The main results of [Str2] and [Str3] respectively classified these posets, thereby classifying the newly expanded types of minuscule representations in the process. The following table summarizes this classification.

---

[*]This paper is a re-envisioning, modernization, and completion of the 1993 manuscript "Poset Partitions and Minuscule Representations: External Constructions of Lie Representations, Part I" by Robert A. Proctor.

|  | Kac–Moody | Positive Borel |
|---|---|---|
| Finite dimensional, highest weight representations | Minuscule posets [Pr1] (Finite type) | $d$-complete posets [Pr2, Ste] (Finite, affine, indefinite types) |
| Infinite dimensional, non-highest weight representations | Full heaps [Gr1] (Affine type) | Filters of full heaps (Affine type) |

The Γ-colored minuscule posets are those occupying the Kac–Moody column and the Γ-colored $d$-complete posets are those occupying the Borel column. The classic minuscule representations of simple Lie algebras are isomorphic to those produced by minuscule posets in the top left cell of this table.

In this paper we work exclusively with finite posets and finite dimensional representations, so any reference to Γ-colored minuscule or Γ-colored $d$-complete posets will be in reference to the first row of the above table. We consider symmetric powers of the various kinds of minuscule representations considered in [Str1]. Let $P$ be a Γ-colored minuscule or Γ-colored $d$-complete poset and let $V$ be its corresponding Kac–Moody or Borel representation. Our main objects of interest in this paper are vectors in $S^m V$ for integers $m \geq 1$. Our main representation theoretic result is Theorem 9.5 where we restrict to simply laced simple Lie algebras $\mathfrak{g}$; here $V$ is isomorphic to one of the classic minuscule representations $V(\lambda)$ for minuscule highest weight $\lambda$. In this result we produce bases for representations isomorphic to the irreducible representations $V(m\lambda)$ of $\mathfrak{g}$ indexed by $m$-multichains of weights in the weight lattice of $V(\lambda)$. C.S. Seshadri [Ses] showed bases indexed in this fashion to exist, and we state a version of that result here for reference.

**Theorem 1.1** (Seshadri). *Let $V(\lambda)$ be a minuscule representation of a simple Lie algebra with highest weight $\lambda$. Then for $m \geq 1$ there exists a weight basis for the irreducible representation of highest weight $m\lambda$ that is indexed in a natural fashion by $m$-multichains in the weight lattice for $V(\lambda)$.*

Moreover, we produce this basis over the integers and in a more general setting. Our constructions and results are entirely combinatorial and we make no algebraic implications until Section 9; no algebraic geometry is needed. For example, our basis vectors are produced from $m$-flags $I_1 \supseteq \cdots \supseteq I_m$ of order ideals of $P$ which are in natural one-to-one correspondence with $m$-multichains in the weight lattice of $V(\lambda)$. As in [Str1], our constructions here concretely describe the actions of Lie algebra generators and are independent of Lie type. (We note other recent work [DDW] in types $E_6$ and $E_7$.) Key definitions to describe our weight bases of "$m$-stackwise vectors" indexed by the $m$-flags noted above are given in Section 3; see Definition 3.3 for our definition of a "standard" poset. Our main combinatorial basis result is Theorem 8.1 where we show that $\Gamma$-colored $d$-complete and $\Gamma$-colored minuscule posets are standard, i.e. that they can be used to produce the integral weight bases of $m$-stackwise vectors indexed by $m$-flags of order ideals.

Though we give no details here, we note that Theorem 8.1 can also be used to provide an important step in the proof that finite $\Gamma$-colored $d$-complete posets are colored hook length posets. These posets correspond to Peterson's $\lambda$-minuscule Weyl group elements [Car] and can be used to produce reduced decompositions of these elements [Pr2, Ste, Str2]. The key insight to show $\Gamma$-colored $d$-complete posets are colored hook length posets is to show that our Borel representations in [Str1] (where $m = 1$) and in this paper are isomorphic to Demazure modules for their corresponding $\lambda$-minuscule Weyl group elements, thus obtaining a description of the Demazure character in terms of $m$-multichains in the Bruhat ideal of the $\lambda$-minuscule Weyl group element corresponding to the poset. The Peterson–Proctor colored hook length identity is [Pr3, Thm. 2] and the key step in its proof to which we refer is Proposition II of that paper; the outline of its proof in [Pr3] relied on descriptions of these Demazure characters given by Lakshmibai [Lak] and Littelmann [Lit] via "admissible weighted $\lambda$-chains." Once it is shown that our representations for $m = 1$ are are isomorphic to Demazure modules, Theorem 8.1 can then be used to replace these references to Lakshmibai and Littelmann, providing a combinatorial proof of Proposition II from [Pr3]. Sections 12 and 13 of [Pr3] give comments about the proof that these posets are (colored) hook length posets and an early history of the notion of "$d$-complete" (prior to this author's contributions).

Both $\Gamma$-colored $d$-complete and $\Gamma$-colored minuscule posets have been used extensively for their algebraic utility. We have already mentioned their usefulness in constructing, characterizing, and classifying representations that possess key properties of minuscule representations of simple Lie algebras, as well as their correspondence with $\lambda$-minuscule Weyl group elements. Uncolored versions of these posets (originally studied by Proctor) have useful applications in $K$-theory; for example [BuSa, IPZ]. We also note the use in [NaOk] of the equivariant $K$-theory of Kac–Moody flag varieties to provide an alternate proof of the Peterson–Proctor hook formula of [Pr3]. The partial order structure of coroots for the corresponding Lie system can be realized with $\Gamma$-colored $d$-complete and $\Gamma$-colored minuscule posets; see Theorem 11 of [Pr1] and Theorem 5.5 of [Ste] and Theorem 9.2 of [Str3]. A survey of many other uses and applications of $d$-complete posets can be found in Sections 1 and 12 of [PrSc] and in the introductions of [Str1, Str2, Str3].

We give definitions in Sections 2 and 3. Section 4 is primarily dedicated to obtaining the linear independence of our $m$-stackwise vectors. The integral spanning proof is unconnected to linear independence and is given in Sections 5–7. The main spanning result is Theorem 6.1; its inductive proof relies on four key lemmas that are proved in Section 7 under the same inductive hypothesis, and so these lemmas should be considered part of the proof of Theorem 6.1. Our main combinatorial result summarizing Sections 4–7 is Theorem 8.1. Section 9 connects to representation theory and concludes in Theorem 9.5 by producing the combinatorially realized basis over the integers for a representation isomorphic to $V(m\lambda)$.

## 2   Basic combinatorial definitions

Let $P$ be a nonempty finite partially ordered set. We use many key terms from [Sta] such as interval, covering relations and the Hasse diagram, order ideal and order filter, distributive lattice, and linear extension. In this paper we consider only connected posets. We also generally follow the notation and terminology used in [Str1]. Letters such as $z, y, x, \ldots$ will denote elements of $P$. If an element $x$ is covered by an element $y$, we denote this by $x \to y$. Two poset elements in a covering relation are *neighbors* in $P$.

Let $\Gamma$ be a finite simple graph. Letters such as $a, b, c, \ldots$ will denote the nodes of $\Gamma$ and we call these nodes *colors*. If there is an edge between distinct colors $a, b \in \Gamma$, then $a$ and $b$ are *adjacent* and we write $a \sim b$. If not, then $a$ and $b$ are *distant*. We will call $\Gamma$ a *Dynkin diagram*.

We *color* the poset $P$ with the colors in $\Gamma$ via a surjective coloring function $\kappa : P \to \Gamma$. A $\Gamma$-*colored poset* is a triple $(P, \Gamma, \kappa)$ which we will shorten to $P$ for brevity. A $\Gamma$-colored poset $P$ may or may not satisfy the following properties.

> (EC) Elements with equal colors are comparable.
>
> (NA) Neighbors have adjacent colors.
>
> (AC) Elements with adjacent colors are comparable.

For $a \in \Gamma$, we define $P_a = \{x \in P \mid \kappa(x) = a\}$. If $x, y \in P_a$ with $x < y$ and there are no elements of color $a$ in the interval $(x, y)$, then $x$ and $y$ are *consecutive elements of color a*.

> (ICE2) For every $a \in \Gamma$, if $x < y$ are consecutive elements of color $a$, then the interval $(x, y)$ contains exactly two elements whose colors are adjacent to $a$.
>
> (UCB1) For every $a \in \Gamma$, if $x$ is maximal in $P_a$, then there is at most one element greater than $x$ with color adjacent to $a$.
>
> (LCB1) For every $a \in \Gamma$, if $x$ is minimal in $P_a$, then there is at most one element less than $x$ with color adjacent to $a$.

Note that ICE2, UCB1, and LCB1 control the total number of elements with adjacent color in either an interval between consecutive elements of the same color or in an upper or lower "frontier" beyond the maximal or minimal element of a given color. We call the latter two properties the *frontier census bounds*.

**Remark 2.1.**   (a)  The assumption that $\Gamma$ is a simple graph makes it a *simply laced* Dynkin diagram. More general Dynkin diagrams and more general coloring axioms can be used; see [Str2] and [Str3].

(b) We presented $\Gamma$ as a graph to emphasize the combinatorial viewpoint, but there is an equivalent way of recognizing its data using integers $\{\theta_{ab}\}_{a,b\in\Gamma}$. For all $a \in \Gamma$ we have $\theta_{aa} = 2$ and for all distinct $a, b \in \Gamma$ we have $\theta_{ab} = \theta_{ba} = -1$ if $a \sim b$ and $\theta_{ab} = \theta_{ba} = 0$ otherwise. Fixing an order to the colors in $\Gamma$ makes $[\theta_{ab}]$ a *generalized Cartan matrix (GCM)* which can be used to define a Kac–Moody algebra $\mathfrak{g}$ and several subalgebras. We do this in Section 9 when we connect our work to representation theory.

Our main colored poset definitions are the following.

**Definition 2.2.** Let $P$ be a $\Gamma$-colored poset. Then $P$ is $\Gamma$-*colored d-complete* if it satisfies EC, NA, AC, ICE2, and UCB1. If it also satisfies LCB1, then it is $\Gamma$-*colored minuscule*.

These colored poset definitions can be given for infinite locally finite posets as well (meaning that every interval is finite). All finite and infinite $\Gamma$-colored $d$-complete posets were classified in [Str2]. All finite and infinite $\Gamma$-colored minuscule posets were classified in [Str3]. These posets can be used to describe certain Weyl group elements as well as roots and coroots of the Lie system associated to $\Gamma$ as in [Gr2, Gr3, Pr1, Pr2, Ste, Str2, Str3], but our main concern will involve their use in creating certain minuscule representations of Lie algebras associated to $\Gamma$ as in [Gr1, Gr3, Str1]; see Section 9.

Let $P$ be a $\Gamma$-colored poset. The set $\mathcal{FI}(P)$ of all filter-ideal *splits* of $P$ is the collection of pairs $(F, I)$ where $F$ is a filter of $P$ and $I = P - F$ is its corresponding ideal. Let $\langle \mathcal{FI}(P) \rangle$ be the $\mathbb{C}$-vector space with basis vectors $\langle F, I \rangle$ corresponding to the splits of $P$. For $a \in \Gamma$ and $(F, I) \in \mathcal{FI}(P)$, we define the *color raising operator* $X_a$ to be $X_a.\langle F, I \rangle = \sum \langle F - \{x\}, I \cup \{x\} \rangle$ and the *color lowering operator* $Y_a$ to be $Y_a.\langle F, I \rangle = \sum \langle F \cup \{x\}, I - \{x\} \rangle$, where these sums are respectively taken over all minimal elements of $F$ and maximal elements of $I$ of color $a$. We extend these operators linearly to $\langle \mathcal{FI}(P) \rangle$. When $P$ satisfies EC, which we always assume, these sums either vanish or are a single term. For brevity, we usually emphasize ideals and write $X_a.\langle F, I \rangle$ as $X_a.\langle I \rangle$ or as $X_a.I$. Throughout most of this paper these operators will be used purely combinatorially. In Section 9 we use these raising and lowering operators to construct representations of Lie algebras.

## 3 Standard posets

Let $m \geq 1$. We use two equivalent notions throughout this paper. The first is a *P-partition bounded by $m$*, which is an order reversing map $\psi : P \to \{0, 1, \ldots, m\}$. The second is an *$m$-flag of order ideals* $I_1 \supseteq I_2 \supseteq \cdots \supseteq I_m$. Given that $\mathcal{FI}(P)$ is a distributive lattice ordered by inclusion of the ideals within the splits, an $m$-flag of order ideals corresponds to an $m$-multichain $(F_1, I_1) \geq (F_2, I_2) \geq \cdots \geq (F_m, I_m)$ in $\mathcal{FI}(P)$. Given a $P$-partition bounded by $m$, an $m$-flag can be constructed where $x \in P$ appears in the first $\psi(x)$ ideals of the $m$-flag. Conversely, given an $m$-flag of order ideals, for $x \in P$ let $\psi(x)$ be the number of ideals containing $x$ in the $m$-flag. Then $\psi : P \to \{0, 1, \ldots, m\}$ is a $P$-partition bounded by $m$. These are corresponding ways of describing the same structure.

Let $\{I_1, \ldots, I_m\}$ be a multiset of $m$ ideals. We call these *$m$-multisets* of ideals and we call $\mathbb{C}$-spans of such multisets *$m$-vectors*. For $a \in \Gamma$ and $k \geq 1$ we define

$$X_a.\{I_1, \ldots, I_m\} = \sum_{j=1}^{m} \{I_1, \ldots, X_a.I_j, \ldots, I_m\} \quad \text{and} \quad \langle a^k \rangle.\{I_1, \ldots, I_m\} = \frac{1}{k!} X_a^k.\{I_1, \ldots, I_m\}.$$

Note that we temporarily fix ideal positions when computing these actions; see Remark 5.1(a).

Let $|P| = p$ and fix a linear extension $x_1 < x_2 < \cdots < x_p$ of $P$ with corresponding colors $a_1, a_2, \ldots, a_p$ for the elements of $P$. Let $I_1 \supseteq I_2 \supseteq \cdots \supseteq I_m$ be an $m$-flag of order ideals with corresponding $P$-partition $\psi : P \to \{0, 1, \ldots, m\}$ bounded by $m$. The *$m$-stackwise operator* for this $P$-partition is $\langle a_p^{\psi(x_p)}, \ldots, a_1^{\psi(x_1)} \rangle$ and the *$m$-stackwise vector* is $\langle a_p^{\psi(x_p)}, \ldots, a_1^{\psi(x_1)} \rangle.\{\emptyset, \ldots, \emptyset\}$, where $\langle a_p^{\psi(x_p)}, \ldots, a_1^{\psi(x_1)} \rangle$ means that we compose right-to-left starting with $\langle a_1^{\psi(x_1)} \rangle$ and ending with $\langle a_p^{\psi(x_p)} \rangle$.

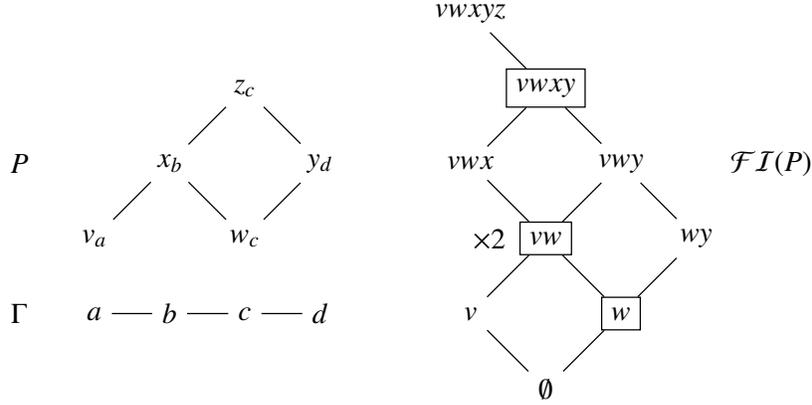

Figure 3.1: A $\Gamma$-colored poset $P$ (left) and its distributive lattice of filter-ideal splits $\mathcal{FI}(P)$ (right). Subscripts for elements in $P$ indicate colors. Splits in $\mathcal{FI}(P)$ are identified by the elements in the ideal in the split without additional symbols. The boxed ideals correspond to the 4-flag $vwxy \supseteq vw \supseteq vw \supseteq w$.

**Example 3.1.** Let $m = 4$ and consider Figure 3.1. As in that figure, our standard practice for notational brevity is to identify splits by the elements the ideal contains without additional symbols. For example, the split $(\{z\}, \{v, w, x, y\})$ is written $vwxy$. The boxed ideals in the figure correspond to the 4-flag $vwxy \supseteq vw \supseteq vw \supseteq w$, and this 4-flag corresponds to the $P$-partition $\psi(v) = 3$, $\psi(w) = 4$, $\psi(x) = 1$, $\psi(y) = 1$, and $\psi(z) = 0$ bounded by 4. Fix the linear extension $v < w < x < y < z$ of $P$. Then the 4-stackwise operator for $\psi$ is $\langle c^0, d^1, b^1, c^4, a^3 \rangle$. Computing the 4-stackwise vector for this operator, we get

$$\langle d^1, b^1, c^4, a^3 \rangle.\{\emptyset, \emptyset, \emptyset, \emptyset\} = \langle d^1, b^1, c^4 \rangle.(4\{v, v, v, \emptyset\}) = \langle d^1, b^1 \rangle.(4\{vw, vw, vw, w\})$$
$$= \langle d^1 \rangle.(12\{vwx, vw, vw, w\}) = 12\{vwxy, vw, vw, w\} + 24\{vwx, vwy, vw, w\} + 12\{vwx, vw, vw, wy\}$$

We note that even though we do not get a single term, the 4-flag we started with appears in this expansion.

**Remark 3.2.** Each $m$-multiset of ideals in the expansion from Example 3.1 contained the same total number of elements with the same *color census* for their elements, i.e. each $m$-multiset contained three elements of color $a$, one of color $b$, four of color $c$, and one of color $d$. This is precisely the same color census of the operator leading to this expansion. This will always be the case in a poset satisfying EC since for all $a \in \Gamma$ the action $X_a.I$ either vanishes or results in a single element of color $a$ added to $I$. So any time a single operator is converted into equivalent linear sums of other operators (notably, as in Section 7 below), each operator in each equivalent expression must have the same color census as the original operator.

Let $m \geq 1$ and let $\mathcal{V}^m$ be the set of all $m$-vectors of the form $\langle b_k^{n_k}, \ldots, b_1^{n_1} \rangle.\{\emptyset, \ldots, \emptyset\}$ for colors $b_1, \ldots, b_k$ and non-negative integers $n_1, \ldots, n_k$. We give the main definition of this paper.

**Definition 3.3.** The poset $P$ is *standard* if for every $m \geq 1$ and every linear extension of $P$ the $m$-stackwise vectors are linearly independent and integrally span $\mathcal{V}^m$.

Theorem 8.1 is the main combinatorial result of this paper and shows that $\Gamma$-colored $d$-complete and $\Gamma$-colored minuscule posets are standard. If $P$ is standard, then the $m$-stackwise vectors form a basis for the vector space spanned by $\mathcal{V}^m$. We show in Section 9 that for $\Gamma$-colored minuscule posets this is a concrete combinatorial construction of a basis shown to exist by C.S. Seshadri, as noted in the Introduction.

## 4 Total orders and linear independence of $m$-stackwise vectors

In this section we continue to work with some fixed $m \geq 1$. We introduce several total orders that we use throughout the rest of the paper. We also prove the linear independence of $m$-stackwise vectors.

The first two total orders we introduce are on ideals and on $m$-multisets of ideals. Fix a linear extension $x_1 < x_2 < \cdots < x_p$ of $P$ and let $I$ and $J$ be ideals of $P$. If $|I| > |J|$ then we define $I < J$ in this total order. If $|I| = |J|$ then we order lexicographically; if $x_{i_1} < x_{i_2} < \cdots < x_{i_k}$ and $x_{j_1} < x_{j_2} < \cdots < x_{j_k}$ are the respective elements of $I$ and $J$ ordered according to the linear extension of $P$, then $I < J$ if $x_{i_\ell} < x_{j_\ell}$ for some $\ell \leq k$ with $x_{i_h} = x_{j_h}$ whenever $h < \ell$. We call this the *gravity order* on ideals; ideals are lower (heavier) in the order when they contain more elements, and they are lower compared to other ideals with the same cardinality when they lead with elements that are lower (heavier) in the linear extension of $P$.

The order on $m$-multisets of ideals is similar. Let $\{I_1, \ldots, I_m\}$ and $\{J_1, \ldots, J_m\}$ be $m$-multisets of ideals of $P$ with $I_1 \leq \cdots \leq I_m$ and $J_1 \leq \cdots \leq J_m$ in the gravity order on ideals. If $\sum_{k=1}^m |I_k| > \sum_{k=1}^m |J_k|$ then we define $\{I_1, \ldots, I_m\} < \{J_1, \ldots, J_m\}$ in this total order. If these cardinality sums are equal then we order lexicographically; we have $\{I_1, \ldots, I_m\} < \{J_1, \ldots, J_m\}$ if $I_\ell < J_\ell$ for some $\ell \leq m$ with $I_h = J_h$ whenever $h < \ell$. We likewise call this the *gravity order* on $m$-multisets of ideals.

**Remark 4.1.** We now always order $m$-multisets lowest to highest according to the gravity order on ideals, and we order $m$-vectors lowest to highest according to the gravity order on $m$-multisets. We ordered the output 4-stackwise vector accordingly in Example 3.1. Note that if $I_1 \supseteq \cdots \supseteq I_m$ is an $m$-flag of ideals, then $I_1 \leq \cdots \leq I_m$ in the gravity order on ideals. Thus the corresponding $m$-multiset of ideals written according to the gravity order is $\{I_1, \ldots, I_m\}$. This is also illustrated in Example 3.1 and the 4-multiset of ideals corresponding to the 4-flag we started with belongs to the leading term; we see in Lemma 4.4 that this is no coincidence.

To prove linear independence of $m$-stackwise vectors we need just two poset coloring axioms. The first is EC and the second is a weakening of NA in which we require neighbors to have different (but not necessarily adjacent) colors. We call this axiom ND.

**Remark 4.2.** If $P$ satisfies EC and ND, then all elements of the same color are comparable but never in a covering relation. Hence $X_a^2.I = 0$ for all $a \in \Gamma$ and ideals $I$. Because of this, whenever $X_a^k.\{I_1, \ldots, I_m\} \neq 0$ for $k > 1$ and $m$-multiset of ideals $\{I_1, \ldots, I_m\}$, every $m$-multiset in each resulting term of this action is a distribution of $k$ elements colored $a$ to $k$ different ideals in $\{I_1, \ldots, I_m\}$.

We call any sequence of colors a *word*. We remark that $b_k, \ldots, b_k, \ldots, b_2, \ldots, b_2, b_1, \ldots, b_1$ (where group $b_i$ appears $n_i$ times for $1 \leq i \leq k$) and $b_k^{n_k}, \ldots, b_2^{n_2}, b_1^{n_1}$ are equal as words, though their corresponding operators $\langle b_k, \ldots, b_k, \ldots, b_2, \ldots, b_2, b_1, \ldots, b_1 \rangle$ and $\langle b_k^{n_k}, \ldots, b_2^{n_2}, b_1^{n_1} \rangle$ differ by a nonzero constant if $n_i > 1$ for any $i$. For notational brevity we will sometimes use bold letters such as $\mathbf{a}$ to denote words. A word corresponding to an $m$-stackwise vector is an *m-stackwise word*. In Section 3 we saw that $m$-flags have corresponding $m$-stackwise vectors and thus have corresponding $m$-stackwise words.

For the remainder of this section it will be useful to consider a particular $m$-multiset of ideals arising from a word's operator acting on $\{\emptyset, \ldots, \emptyset\}$. For an $m$-multiset $\{I_1, \ldots, I_m\}$ listed by the gravity order and $a \in \Gamma$ we define $X_a^{\min}.\{I_1, \ldots, I_m\} = \{I_1, \ldots, X_a.I_g, \ldots, I_m\}$ where $g$ is minimal such that $X_a.I_g \neq 0$, or $0$ otherwise. We say a word $b_k, \ldots, b_2, b_1$ *grows well for* $m$ if $X_{b_k}^{\min} \cdots X_{b_2}^{\min} X_{b_1}^{\min}.\{\emptyset, \ldots, \emptyset\} \neq 0$.

**Lemma 4.3.** *Suppose P satisfies EC and ND. Let $I_1 \supseteq \cdots \supseteq I_m$ be an m-flag with word $\mathbf{a} = a_p^{n_p}, \ldots, a_2^{n_2}, a_1^{n_1}$.*
  *(a) Let $1 \leq j \leq p$ and suppose that $(X_{a_{j-1}}^{\min})^{n_{j-1}} \cdots (X_{a_2}^{\min})^{n_2} (X_{a_1}^{\min})^{n_1}.\{\emptyset, \ldots, \emptyset\} = \{J_1, \ldots, J_m\} \neq 0$. Then we have $(X_{a_j}^{\min})^{n_j}.\{J_1, \ldots, J_m\} = \{J_1 \cup \{x_j\}, \ldots, J_{n_j} \cup \{x_j\}, J_{n_j+1}, \ldots, J_m\}$.*
  *(b) We have $(X_{a_p}^{\min})^{n_p} \cdots (X_{a_2}^{\min})^{n_2} (X_{a_1}^{\min})^{n_1}.\{\emptyset, \ldots, \emptyset\} = \{I_1, \ldots, I_m\}$ and so $\mathbf{a}$ grows well for m.*

*Proof.* We prove (a) by induction on $j$. For $j = 1$ the nonzero starting term is $\{\emptyset, \ldots, \emptyset\}$. Since $a_1$ is the color of the unique minimal element $x_1$ of $P$ of color $a_1$ and $n_1 \leq m$, by Remark 4.2 we have $(X_{a_1}^{\min})^{n_1}.\{\emptyset, \ldots, \emptyset\} = \{x_1, \ldots, x_1, \emptyset, \ldots, \emptyset\}$. Here there are $n_1$ copies of $x_1$ in this $m$-multiset of ideals and $m - n_1$ copies of $\emptyset$. (Note the possibility that $n_1 = 0$, which does not affect this argument.) This proves the base case.

Now suppose $1 < j \leq p$ and that the result holds for every $h < j$. Let $\{J_1, \ldots, J_m\}$ be the $m$-multiset resulting from step $j - 1$. Note that $J_1 \supseteq \cdots \supseteq J_m$ and so Remark 4.1 shows $\{J_1, \ldots, J_m\}$ is listed in the gravity order. For every $h < j$ with $x_h < x_j$ in $P$, we see that $x_h$ was added to the first $n_h$ ideals in the final step of $(X_{a_h}^{\min})^{n_h} \cdots (X_{a_1}^{\min})^{n_1}.\{\emptyset, \ldots, \emptyset\}$. Moreover, the element $x_j$ has not be added to any ideal in any of the first $j - 1$ steps. So $x_j$ may be added to the first $\min\{n_h \mid h < j$ with $x_h < x_j$ in $P\}$ ideals of $\{J_1, \ldots, J_m\}$. The $P$-partition bounded by $m$ ensures that this number is at least $n_j$ and at most $m$. Hence by Remark 4.2 we see that $(X_{a_j}^{\min})^{n_j}.\{J_1, \ldots, J_m\} = \{J_1 \cup \{x_j\}, \ldots, J_{n_j} \cup \{x_j\}, J_{n_j+1}, \ldots, J_m\}$, so (a) is proved.

Applying Part (a) for $1 \leq j \leq p$ shows that the final action $(X_{a_p}^{\min})^{n_p} \cdots (X_{a_2}^{\min})^{n_2} (X_{a_1}^{\min})^{n_1}.\{\emptyset, \ldots, \emptyset\}$ will consist of an $m$-multiset of ideals where the first $n_1$ ideals contain $x_1$, the first $n_2$ ideals contain $x_2$, etc, and finally, the first $n_p$ ideals contain $x_p$. This is the precise description of $I_1 \supseteq \cdots \supseteq I_m$ that led to the $m$-stackwise word $\mathbf{a}$, and so (b) holds. □

The particular term arising from Lemma 4.3 gives the leading $m$-multiset of ideals in a stackwise action.

**Lemma 4.4.** *Suppose P satisfies EC and ND. Let $I_1 \supseteq \cdots \supseteq I_m$ be an m-flag with word $\mathbf{a} = a_p^{n_p}, \ldots, a_2^{n_2}, a_1^{n_1}$. Then the m-multiset of ideals in the leading term of $\langle a_p^{n_p}, \ldots, a_2^{n_2}, a_1^{n_1} \rangle.\{\emptyset, \ldots, \emptyset\}$ is $\{I_1, \ldots, I_m\}$.*

*Proof.* Throughout this proof we ignore coefficients since all coefficients would be positive at each step (i.e. no terms cancel) and this result only concerns the $m$-multisets appearing in successive steps of the action.

We consider the action $\langle a_p^{n_p}, \ldots, a_2^{n_2}, a_1^{n_1} \rangle.\{\emptyset, \ldots, \emptyset\}$ one color group at a time. Since $\mathbf{a}$ grows well for $m$ by Lemma 4.3(b), we note that for every $1 \leq j \leq p$ the $m$-multiset $(X_{a_j}^{\min})^{n_j} \cdots (X_{a_2}^{\min})^{n_2} (X_{a_1}^{\min})^{n_1}.\{\emptyset, \ldots, \emptyset\}$ appears in the expansion of $\langle a_j^{n_j}, \ldots, a_2^{n_2}, a_1^{n_1} \rangle.\{\emptyset, \ldots, \emptyset\}$. If the result after each color group is a single $m$-multiset of ideals, then Lemma 4.3(b) also shows it must be $\{I_1, \ldots, I_m\}$, so the result holds.

Now suppose there is a step of the action where more than one $m$-multiset of ideals arises. Assume this happens first at color group $a_h^{n_h}$ so that $\langle a_{h-1}^{n_{h-1}}, \ldots, a_1^{n_1} \rangle.\{\emptyset, \ldots, \emptyset\}$ is a single term. Since $\langle a_1^{n_1} \rangle.\{\emptyset, \ldots, \emptyset\}$ is a single term, we see that $h > 1$. Let $\{J_1, \ldots, J_m\}$ be the $m$-multiset of ideals in this term at step $h-1$, which must be the $m$-multiset described by Lemma 4.3(a). We note that $J_1 \supseteq \cdots \supseteq J_m$. By assumption $\langle a_h^{n_h} \rangle.\{J_1, \ldots, J_m\}$ results in two or more $m$-multisets of ideals. By Lemma 4.3(a) one resulting $m$-multiset is $\{J_1 \cup \{x_h\}, \ldots, J_{n_h} \cup \{x_h\}, J_{n_h+1}, \ldots, J_m\}$, which is written according to the gravity order. Every other resulting $m$-multiset must result from a distribution of $n_h$ elements to different ideals in $\{J_1, \ldots, J_k\}$ by Remark 4.2. Let $\{K_1, \ldots, K_m\}$ be another $m$-multiset arising in the expansion of $\langle a_h^{n_h} \rangle.\{J_1, \ldots, J_m\}$ and note both it and $\{J_1 \cup \{x_h\}, \ldots, J_{n_h} \cup \{x_h\}, J_{n_h+1}, \ldots, J_m\}$ have $n_1 + \cdots + n_h$ elements. Suppose that $\{K_1, \ldots, K_m\}$ differs from $\{J_1 \cup \{x_h\}, \ldots, J_{n_h} \cup \{x_h\}, J_{n_h+1}, \ldots, J_m\}$ first at index $r$ and note that $1 \leq r \leq n_h$ since these two $m$-multisets would be identical if they agree in the first $n_h$ ideals. Hence we have $K_\ell = J_\ell \cup \{x_h\}$ whenever $\ell < r$. Since $\{K_1, \ldots, K_m\}$ arose from acting on $\{J_1, \ldots, J_m\}$ with $\langle a_h^{n_h} \rangle$, we see that $K_r = J_s$ for some $s \geq r$ or $K_r = J_s \cup \{y\}$ for some $s \geq r$ and element $y$ of color $a_h$. If $K_r = J_s$ for some $s \geq r$, then since $J_1 \supseteq \cdots \supseteq J_m$ we have $J_r \cup \{x_h\} \supseteq J_r \supseteq J_s$. Hence $J_r \cup \{x_h\} < K_r$ in the gravity order by cardinality. Now suppose that $K_r = J_s \cup \{y\}$ for some $s \geq r$ and element $y$ of color $a_h$. If $J_s = J_r$, then we note by EC that $y = x_h$ and so $K_r = J_r \cup \{x_h\}$. This contradicts that $\{J_1 \cup \{x_h\}, \ldots, J_{n_h} \cup \{x_h\}, J_{n_h+1}, \ldots, J_m\}$ and $\{K_1, \ldots, K_m\}$ differ at the $r^{\text{th}}$ ideal. Hence since $J_1 \supseteq \cdots \supseteq J_m$ we see that $J_r$ properly contains $J_s$. Thus $J_r \cup \{x_h\} < K_r$ in the gravity order, again by cardinality. In both cases $\{J_1 \cup \{x_h\}, \ldots, J_{n_h} \cup \{x_h\}, J_{n_h+1}, \ldots, J_m\} < \{K_1, \ldots, K_m\}$ in the gravity order, with cardinality being the distinguishing factor in the lexicographic comparison of ideals within each $m$-multiset.

Sequential actions for $j > h$ on each resulting higher $m$-multiset in the gravity order can never lower them to the Lemma 4.3(a) term. This is because all $n_j$ elements are added to the first $n_j$ ideals in the Lemma 4.3(a) term, always resulting in the same or possibly more separated cardinality comparison with other terms. Moreover, all of these terms always have the same initial ideals whenever beginning cardinalities match since these ideals grow by a unique element of color $a_j$, namely $x_j$. Thus the final Lemma 4.3(a) term is the $m$-multiset of ideals in the leading term of $\langle a_p^{n_p}, \ldots, a_1^{n_1} \rangle.\{\emptyset, \ldots, \emptyset\}$. From Lemma 4.3(b) we see this $m$-multiset of ideals is $\{I_1, \ldots, I_m\}$. $\square$

**Remark 4.5.** Suppose $P$ satisfies EC and ND. Remark 4.2 and Lemmas 4.3 and 4.4 combine to describe where the distribution of elements in any $m$-multiset of ideals must occur in an $m$-stackwise vector $\langle a_p^{n_p}, \ldots, a_1^{n_1} \rangle.\{\emptyset, \ldots, \emptyset\}$ coming from an $m$-flag $I_1 \supseteq \cdots \supseteq I_m$. A color group operator $\langle a_j^{n_j} \rangle$ for $1 \leq j \leq p$ can only add one element at a time to any ideal in any $m$-multiset appearing at that stage of the operation, and if it does add an element, this element has color $a_j$. The leading term receives $n_j$ copies of $x_j$ and this is the greatest element of color $a_j$ able to be added from $P$ to any ideal at this point since elements must be added from the bottom of $P$. In other words, if this is the $k^{\text{th}}$ time that color $a_j$ has appeared in a group for the operator $\langle a_j^{n_j}, \ldots, a_1^{n_1} \rangle.\{\emptyset, \ldots, \emptyset\}$, then only the lowest $k$ elements of color $a_j$ in $P$ could have possibly been added to any ideal at this point in the operation. Thus all ideals appearing in all $m$-multisets of ideals in the expansion of $\langle a_p^{n_p}, \ldots, a_1^{n_1} \rangle.\{\emptyset, \ldots, \emptyset\}$ are contained in $I_1$.

We can now obtain linear independence, which only requires that $P$ satisfies EC and ND. Before we do,

we comment on a converse. Let $m \geq 2$. Suppose that $P$ does not satisfy EC. Let $x, y \in P$ be incomparable having the same color $b$. Let $F$ be the principal filter generated by $x$ and $y$ and let $I$ be its corresponding ideal. Fix a linear extension $x_1 < \cdots < x_j < x < y < \cdots < x_p$ of $P$, where $I = \{x_1, \ldots, x_j\}$. Let $\psi_1 : P \to \{0, 1, \ldots, m\}$ be the $P$-partition bounded by $m$ where $\psi_1(u) = m$ if $u \in I$ and $\psi_1(x) = 1$ and $\psi_1(y) = 1$ and $\psi(v) = 0$ if $v \in F - \{x, y\}$. Let $\psi_2 : P \to \{0, 1, \ldots, m\}$ be equal to $\psi_1$ on all elements except $x$ and $y$ and define $\psi_2(x) = 2$ and $\psi_2(y) = 0$. Then $\psi_1$ produces the $m$-stackwise vector $\langle b, b, a_j^m, \ldots, a_1^m \rangle.\{\emptyset, \ldots, \emptyset\}$ and $\psi_2$ produces the $m$-stackwise vector $\langle b^2, a_j^m, \ldots, a_1^m \rangle.\{\emptyset, \ldots, \emptyset\}$, so the $m$-stackwise vectors are linearly dependent. If $P$ fails ND with neighbors $x \to y$ of color $b$, then a similar construction using the principal filter generated by $x$ generates the same linearly dependent scenario. Hence EC and ND fully characterizes linear independence.

**Proposition 4.6.** *Let $P$ be a $\Gamma$-colored poset. The $m$-stackwise vectors are linearly independent for every $m \geq 1$ and every linear extension of $P$ if and only if $P$ satisfies EC and ND.*

*Proof.* The remarks preceding the proposition statement show that EC and ND are necessary for the linear independence of $m$-stackwise vectors for every $m \geq 1$ and every linear extension of $P$. Now assume $P$ satisfies EC and ND. Let $m \geq 1$ and fix a linear extension of $P$. By Lemma 4.4 all $m$-stackwise vectors have distinct leading terms since they originate from distinct $m$-flags that make up these leading terms. Thus these vectors can be arranged in upper triangular form with respect to their expansions into $m$-multisets of ideals in the gravity order, so they are linearly independent. □

We conclude this section by introducing a total order on words. We induct on this order in the proof of Theorem 6.1. We first fix an expansion of the linear extension of $P$ to a total order on $P \cup \Gamma$ where every element of $P$ is greater than every element of $\Gamma$. If $\mathbf{a} = a_k, \ldots, a_2, a_1$ is a word of length $k$ and $1 \leq h \leq k$, then we define $R_h(\mathbf{a})$ to be $a_h, \ldots, a_2, a_1$ and call it an *initial segment* of $\mathbf{a}$ from the right.

**Definition 4.7.** Let $\mathbf{a} = a_j, \ldots, a_2, a_1$ and $\mathbf{b} = b_k, \ldots, b_2, b_1$ be words. If $j < k$ then we define $\mathbf{a} < \mathbf{b}$. Now suppose $j = k$. Let $\ell, h \leq k$ be such that $R_\ell(\mathbf{a})$ and $R_h(\mathbf{b})$ are the longest initial segments of these respective words that grow well for $m$. Replace these initial segments of colors with the corresponding elements of $P$ that were added to ideals in the leading term as these segments grew well; call these new entities $[\mathbf{a}]$ and $[\mathbf{b}]$, respectively. Then we define $\mathbf{a} < \mathbf{b}$ if $[\mathbf{a}] < [\mathbf{b}]$ lexicographically from the right within the total order on $P \cup \Gamma$. That is, if $[\mathbf{a}]_r$ is the $r^{\text{th}}$ position of $[\mathbf{a}]$ from the right, then $[\mathbf{a}] < [\mathbf{b}]$ if $[\mathbf{a}]_\ell < [\mathbf{b}]_\ell$ for some $\ell \leq j$ with $[\mathbf{a}]_h = [\mathbf{b}]_h$ whenever $h < \ell$. This is the *total word order* for $m$ subject to the total order on $P \cup \Gamma$.

We illustrate the total word order with an example.

**Example 4.8.** Let $m = 2$ and refer again to Figure 3.1. We again use the linear extension $v < w < x < y < z$ of $P$ from Example 3.1. Consider the word $\mathbf{a} = d^1, b^1, c^2, a^2$ corresponding to the $P$-partition $\psi(v) = 2$, $\psi(w) = 2, \psi(x) = 1, \psi(y) = 1$, and $\psi(z) = 0$ bounded by 2. This word grows well for 2 by Lemma 4.3(b) so the replaced initial segment as in Definition 4.7 is $[\mathbf{a}] = y, x, w, w, v, v$.

Now consider the word $\mathbf{b} = d, c, d, b, c, a$ which has the same length as $\mathbf{a}$. Note that

$$\langle d,c,d,b,c,a\rangle.\{\emptyset,\emptyset\} = \langle d,c,d,b,c\rangle.(2\{v,\emptyset\}) = \langle d,c,d,b\rangle.(2\{vw,\emptyset\} + 2\{v,w\}) = \langle d,c,d\rangle.(2\{vwx,\emptyset\})$$
$$= \langle d,c\rangle.(2\{vwxy,\emptyset\}) = \langle d\rangle.(2\{vwxyz,\emptyset\} + 2\{vwxy,w\}) = 2\{vwxy,wy\}$$

The leading term does not vanish until the last action, and so the longest initial portion of this word that grows well for 2 is $c, d, b, c, a$. Hence the replaced initial segment as in Definition 4.7 is $[\mathbf{b}] = d, z, y, x, w, v$. Now we compare lexicographically and note that $\mathbf{a} < \mathbf{b}$ since $[\mathbf{a}]_1 = v = [\mathbf{b}]_1$ in the first term from the right but $[\mathbf{a}]_2 = v < w = [\mathbf{b}]_2$ in the second term.

## 5 Operator identities

Since $\Gamma$-colored $d$-complete and $\Gamma$-colored minuscule posets satisfy EC and ND, we showed that their $m$-stackwise vectors are linearly independent in Proposition 4.6. We now focus on showing these vectors integrally span $\mathcal{V}^m$. This will be achieved in Theorem 8.1 by applying Theorem 6.1. Our first step toward that goal is to prove three operator identities for the raising action on $m$-multisets of ideals, Propositions 5.2, 5.4, and 5.6. In this section we consider this action on ordered $m$-tuples of ideals using the same definitions

$$X_a.(I_1,\ldots,I_m) = \sum_{j=1}^{m}(I_1,\ldots,X_a.I_j,\ldots,I_m) \quad \text{and} \quad \langle a^k\rangle.(I_1,\ldots,I_m) = \frac{1}{k!}X_a^k.(I_1,\ldots,I_m)$$

for $a \in \Gamma$ and $k \in \mathbb{N}$ used for these actions on $m$-multisets of ideals. In other words, in this section we consider ideal positions as fixed and do not reorder ideals according to the gravity order. We call a position in an $m$-tuple of ideals a *slot* and say that a slot *accepts* $a \in \Gamma$ if the ideal $I$ in that slot satisfies $X_a.I \neq 0$.

We note that while only axioms EC and ND were needed in Section 4 to prove linear independence of $m$-stackwise vectors, we will need EC, NA, AC, and ICE2 for Propositions 5.4 and 5.6. The pair EC and NA produce the effect described for EC and ND in Remark 4.2, while NA is required to apply Proposition 5.2 in the proofs of Lemma 5.3 and Proposition 5.6. The axiom AC ensures order matters when augmenting elements with adjacent colors to an ideal. Finally, ICE2 prevents three element interval chains in $P$ whose top and bottom elements have the same color.

**Remark 5.1.** (a) The operator identities proved in this section for ordered $m$-tuples of ideals still hold for $m$-multisets of ideals that have been reordered according to the gravity order. If two sides of an operator identity (such as in Propositions 5.2, 5.4, or 5.6) produce equal outputs for ordered $m$-tuples, then these outputs will still be equal when the orders are rearranged into $m$-multisets according to the gravity order, albeit with perhaps larger coefficients as some $m$-tuples may result in the same $m$-multiset when put in the gravity order. Parts (b) and (c) of this remark describe these coefficients.

(b) Consider an action of an operator $\langle b^k\rangle$ on an $m$-tuple $(I_1,\ldots,I_m)$ of ideals. Suppose that the $m$-tuple of ideals $(I_1,\ldots,I_m)$ contains $k \geq j$ slots that accept the color $b$; the action vanishes if $k < j$. An $m$-tuple $(J_1,\ldots,J_m)$ appears in the expansion $\langle b^j\rangle.(I_1,\ldots,I_m)$ if $J_\ell$ is $I_\ell$ augmented with a single element of color $b$ for $j$ slots and where $J_\ell = I_\ell$ for the remaining $m - j$ slots. Before dividing, such an $m$-tuple appears exactly $j!$ times given by the number of ways of ordering these augmentations, so dividing by $j!$ leaves a coefficient of 1 for $(J_1, \ldots, J_m)$ in the expansion. This observation is used in the proofs in this section.

(c) Coefficients can still be quickly computed when the action is instead executed on the $m$-multiset $\{I_1, \ldots, I_m\}$. Suppose that $k_1, \ldots, k_s$ are the multiplicities of distinct ideals in the $m$-multiset $\{I_1, \ldots, I_m\}$ so that $k_1 + \cdots + k_s = m$, and suppose that each ideal has been augmented $j_1 \leq k_1, \ldots, j_s \leq k_s$ times to produce $\{J_1, \ldots, J_m\}$ so that $j_1 + \cdots + j_s = j$. Then $\{J_1, \ldots, J_m\}$ appears with coefficient $\binom{k_1}{j_1} \cdots \binom{k_s}{j_s}$ in the expansion of $\langle b^j \rangle.\{I_1, \ldots, I_m\}$ into $m$-multisets of ideals.

**Proposition 5.2.** *Let P be a $\Gamma$-colored poset and let $m \geq 1$ and $p, q \geq 0$. Suppose that $b, c \in \Gamma$ are never the colors of neighbors in P. Then the following operator identity holds on $m$-tuples and $m$-multisets of ideals.*

$$\langle c^q, b^p \rangle = \langle b^p, c^q \rangle$$

Note that the hypotheses are met if $P$ satisfies NA and $b$ and $c$ are distant in $\Gamma$.

*Proof.* Since $b$ and $c$ are never the colors of neighbors in $P$, a slot accepting an element of color $b$ or $c$ never affects whether that slot accepts an element of the other color. So the result does not depend on whether we attempt to distribute $p$ elements of color $b$ and then $q$ elements of color $c$, or vice versa. □

A quick calculation for $c \in \Gamma$ and $p, q \geq 0$ shows $\langle c^q, c^p \rangle = \binom{p+q}{p} \langle c^{p+q} \rangle$. We call this the *gluing identity* with the *gluing coefficient* $\binom{p+q}{p}$. It applies whether the action is on $m$-tuples or $m$-multisets of ideals.

**Lemma 5.3.** *Let P satisfy EC, NA, AC, and ICE2, and let $m \geq 1$ and $q > p \geq 0$. Then for distinct colors $b, c \in \Gamma$ the following operator identity holds on $m$-tuples and $m$-multisets of ideals.*

$$\sum_{k=0}^{q} (-1)^k \langle c^{q-k}, b^p, c^k \rangle = 0$$

*Proof.* Our proof is for $m$-tuples of ideals; the result also holds for $m$-multisets of ideals by Remark 5.1(a).

First suppose that $b$ and $c$ are distant. Then applying Proposition 5.2 and gluing identities, the identity we must prove becomes

$$\sum_{k=0}^{q} (-1)^k \binom{q}{k} \langle c^q, b^p \rangle = 0$$

The operator $\langle c^q, b^p \rangle$ does not depend on $k$ and the coefficients are alternating binomial coefficients summing to zero, so the identity holds when $b$ and $c$ are distant.

Now suppose $b$ and $c$ are adjacent. Consider the action of the left hand side on an $m$-tuple of ideals $\mathcal{M} = (I_1, \ldots, I_m)$. We must distribute $q$ total elements colored $c$ and $p$ total elements colored $b$. Fix some $0 \leq k \leq q$ and suppose $\mathcal{M}'$ is an $m$-tuple arising from the $k^{\text{th}}$ term in the sum. We count how many total ways $\mathcal{M}'$ can arise.

Notice that two elements of color $c$ without an element of color $b$ cannot be added to a single slot since $P$ satisfies EC and neighbors have distinct colors. Also, no progression of an element colored $c$ followed by an element colored $b$ followed by an element colored $c$ can be added to a single slot. This is because elements of color $c$ form a chain by EC and consecutive elements of color $c$ must have at least two elements between them by ICE2. Hence if an ideal in $\mathcal{M}'$ differs from its corresponding slot ideal in $\mathcal{M}$, then it must be from either: (i) adding a single element of color $c$, (ii) adding a single element of color $b$, (iii) adding an element of color $c$ before an element of color $b$, or (iv) adding an element of color $b$ before an element of color $c$. Since $b$ and $c$ are adjacent and $P$ satisfies AC, groups (iii) and (iv) are distinct, non-overlapping possibilities.

We note that all distribution choices to produce $\mathcal{M}'$ occur at the first operator $\langle c^k \rangle$. Once those elements have been distributed, the $p$ elements of color $b$ must go into their predetermined slots, and the remaining $q - k$ elements of color $c$ must go into the slots that have not met their quota of incoming elements of color $c$. Among the groups identified above, groups (i) and (iii) involve elements potentially being added by the operator $\langle c^k \rangle$. Let $s$ be the total number of elements in group (i) and let $t$ be the total number of elements colored $c$ in group (iii). These numbers $s$ and $t$ only depend on $\mathcal{M}'$, not on $k$. However, we note that $t \leq k$ since all $t$ elements must be distributed by $\langle c^k \rangle$. Note also that $k - t \leq s$ since $k - t$ represents the number of slots that must be filled by the operator $\langle c^k \rangle$ that do not have an element of color $b$ to arrive from $\langle b^p \rangle$, while some of the $s$ elements may be added by $\langle c^{q-k} \rangle$. Observe that $s + t \leq q$ since there are $q$ total elements of color $c$ added to ideals and $s + t$ may not represent all of them. So the summation we are considering can only produce $\mathcal{M}'$ if $t \leq k \leq s + t$.

The $t$ elements colored $c$ in (iii) must be added in the first step $\langle c^k \rangle$, which leaves $k - t$ remaining slots for $\langle c^k \rangle$ to fill. Of these, a total of $s$ slots must eventually wind up with only an element of color $c$, so there are $\binom{s}{k-t}$ ways to fill them in this first step. As noted in Remark 5.1(b), each of these $\binom{s}{k-t}$ preliminary $m$-tuples will have coefficient 1 after the operation $\langle c^k \rangle$.

Once $\langle c^k \rangle$ is complete, everything is determined; the $p$ elements colored $b$ fill their slots and the remaining elements colored $c$ fill their slots, each with coefficient 1 by Remark 5.1(b). So the coefficient on $\mathcal{M}'$ at step $k$ is $(-1)^k \binom{s}{k-t}$. Since $s$ and $t$ are independent of $k$ and $t \leq k \leq s + t$, the total coefficient on $\mathcal{M}'$ is

$$\sum_{k=t}^{s+t} (-1)^k \binom{s}{k-t} = \sum_{h=0}^{s} (-1)^{h+t} \binom{s}{h}$$

Since $q$ elements of color $c$ and $p$ elements of color $b$ are added to slots and $q > p$, we must have $s > 0$. Hence this is a vanishing alternating sum of binomial coefficients. $\square$

**Proposition 5.4.** *Let $P$ satisfy EC, NA, AC, and ICE2, and let $m \geq 1$ and $q > p \geq 0$. Then for distinct colors $b, c \in \Gamma$ the following operator identity holds on m-tuples and m-multisets of ideals.*

$$\langle c^q, b^p \rangle = \sum_{k=0}^{q-1} (-1)^k \langle c^{q-1-k}, b^p, c^{k+1} \rangle$$

*Proof.* Re-indexing from $k = 1$ to $q$ and subtracting the right hand side to the left, this identity becomes

$$0 = \langle c^q, b^p \rangle + \sum_{k=1}^{q} (-1)^k \langle c^{q-k}, b^p, c^k \rangle = \sum_{k=0}^{q} (-1)^k \langle c^{q-k}, b^p, c^k \rangle$$

This is the identity proved in Lemma 5.3. □

We provide combinatorial identities needed for the final operator identity in this section.

**Lemma 5.5.** *Suppose $p, q, r \geq 1$ and $r - p \geq 0$ and $1 \leq k \leq q$.*

*(a) We have $\binom{q + r - p}{q}(1 + q + r - p) = (q + 1)\binom{(q+1) + r - p}{q+1}$.*

*(b) We have $\binom{q + r - p}{q - k}(1 + q + r - p + k) + \binom{q + r - p}{q - k + 1}k = (q + 1)\binom{(q+1) + r - p}{(q+1) - k}$.*

The restrictions on $p, q, r,$ and $k$ ensure these binomial coefficients are always defined. Note that (a) is a special case of (b) with $k = 0$; we list it separately to ensure $\binom{q+r-p}{q-k+1}$ is defined.

*Proof.* For (a) we have

$$\frac{(q + r - p)!}{q!(r - p)!}((q + 1) + r - p) = \frac{((q + 1) + r - p)!}{q!(r - p)!} = (q + 1)\frac{((q + 1) + r - p)!}{(q + 1)!(r - p)!}$$

The left and right hand sides are the respective left and right hand sides in the stated equation for (a).

For (b) we have

$$\frac{(q + r - p)!}{(q - k)!(k + r - p)!}((q + 1) + r - p + k) + \frac{(q + r - p)!}{((q + 1) - k)!(k + r - p - 1)!}k$$

$$= \frac{((q + 1) + r - p)!}{(q - k)!(k + r - p)!} + \frac{(q + r - p)!}{(q - k)!(k + r - p)!}k + \frac{(q + r - p)!}{((q + 1) - k)!(k + r - p - 1)!}k$$

$$= ((q + 1) - k)\frac{((q + 1) + r - p)!}{((q + 1) - k)!(k + r - p)!}$$

$$+ \left[((q + 1) - k)\frac{(q + r - p)!}{((q + 1) - k)!(k + r - p)!} + \frac{(q + r - p)!}{((q + 1) - k)!(k + r - p - 1)!}\right]k$$

$$= ((q + 1) - k)\frac{((q + 1) + r - p)!}{((q + 1) - k)!(k + r - p)!}$$

$$+ \left[\frac{(q + r - p)!}{((q + 1) - k)!(k + r - p)!}((q + 1) - k + (k + r - p))\right]k$$

$$= ((q + 1) - k)\frac{((q + 1) + r - p)!}{((q + 1) - k)!(k + r - p)!} + \frac{((q + 1) + r - p)!}{((q + 1) - k)!(k + r - p)!}k$$

$$= (q + 1)\frac{((q + 1) + r - p)!}{((q + 1) - k)!(k + r - p)!}$$

The first and last parts are the respective left and right hand sides in the stated equation. □

**Proposition 5.6.** *Let P satisfy EC, NA, AC, and ICE2, and let $m \geq 1$ and $p, q, r \geq 0$ and $r - p \geq 0$. Then for distinct colors $b, c \in \Gamma$ the following operator identity holds on m-tuples and m-multisets of ideals.*

$$\langle c^q, b^p, c^r \rangle = \sum_{k=0}^{\min\{p,q\}} \binom{q+r-p}{q-k} \langle b^{p-k}, c^{q+r}, b^k \rangle$$

*Proof.* Our proof is for *m*-tuples of ideals; the result also holds for *m*-multisets of ideals by Remark 5.1(a).

This identity is a tautology if $q = 0$ and is the gluing identity if $p = 0$, so assume that $p, q > 0$. Since $r - p \geq 0$ and $0 \leq k \leq \min\{p, q\}$, we have $0 \leq q - k \leq q + r - p$. Hence the binomial coefficient in the summation identity is always defined.

First suppose that $b$ and $c$ are distant. Proposition 5.2 then applies, so changing the order of colors and applying gluing coefficients, the identity we must prove is

$$\binom{q+r}{q} = \sum_{k=0}^{\min\{p,q\}} \binom{q+r-p}{q-k}\binom{p}{k}$$

Suppose we have $q + r$ objects and want to choose $q$ of them. This can be done directly, resulting in $\binom{q+r}{q}$ possibilities. Alternatively, split these objects into two piles, one of size $p$ and one of size $q + r - p$. We choose $0 \leq k \leq p$ objects from the first pile in $\binom{p}{k}$ ways and the remaining $q - k$ objects from the second pile in $\binom{q+r-p}{q-k}$ ways for a total of $\binom{q+r-p}{q-k}\binom{p}{k}$ possibilities. Note that $k \leq q$ since we are choosing $q$ objects total, so $0 \leq k \leq \min\{p, q\}$. Summing over all possible $k$ gives the quantity on the right-hand side.

Now suppose that $b$ and $c$ are adjacent. Fix $r, p > 0$ with $r - p \geq 0$. We induct on $q \geq 1$. The base case for $q = 1$ is the identity

$$\langle c, b^p, c^r \rangle = \sum_{k=0}^{1} \binom{1+r-p}{1-k} \langle b^{p-k}, c^{1+r}, b^k \rangle = (1 + r - p)\langle b^p, c^{r+1} \rangle + \langle b^{p-1}, c^{r+1}, b \rangle$$

Suppose that both the left and right sides of this proposed identity act on an *m*-tuple of ideals $\mathcal{M}$. Let $A$ be the resulting expansion for the left hand side and let $B$ and $C$ be the resulting expansions for the two terms occurring on the right hand side, respectively. Suppose $\mathcal{M}'$ occurs in the expansion on either the left or right hand side. We count the number of ways $\mathcal{M}'$ occurs in $A$, $B$, and $C$.

As in the proof of Lemma 5.3, no three element interval chains with colors *c-b-c* or *b-c-b* occur because $P$ satisfies ICE2. So no slot in $\mathcal{M}$ can receive a *c*-then-*b*-then-*c* or a *b*-then-*c*-then-*b* addition.

First suppose $\mathcal{M}'$ has a slot that received an element of color $b$ followed by an element of color $c$. Since $b$ and $c$ are adjacent, by AC we know that action order matters. Hence $\mathcal{M}'$ does not occur in $B$ since slots in $B$ must come from *c*-then-*b* additions. There can be only one such slot that received *b*-then-*c* in $\mathcal{M}'$ since the left hand side ends with a single $\langle c \rangle$ operator and the right hand side begins with a single $\langle b \rangle$ operator. Since $\mathcal{M}$ has an ideal that accepts *b*-then-*c* to produce $\mathcal{M}'$, this ideal was produced directly in $C$ by first adding an element of color $b$ to the required slot, then all remaining $r + 1$ elements colored $c$, then the remaining $p - 1$ elements colored $b$. All choices were predetermined in this sequence of actions, so $\mathcal{M}'$ has coefficient 1 in $C$ by Remark 5.1(b). Similarly, $\mathcal{M}'$ was produced in $A$ directly by first placing $r$ elements of color $c$ into the $r$ slots not requiring $b$ first, then placing the $p$ elements of color of $b$ (including one element into the required $b$-then-$c$ slot), and then the final element of color $c$ into the required slot. Again all choices were predetermined, so $\mathcal{M}'$ also has coefficient 1 in $A$ by Remark 5.1(b).

Now assume that $\mathcal{M}'$ has no slot that received an element of color $b$ followed by an element of color $c$. Let $s$ be the number of slots in $\mathcal{M}'$ that received a $c$-then-$b$ addition. To produce $\mathcal{M}'$ in $A$, the $s$ slots to receive $c$-then-$b$ must be filled in the first action step to ensure they are prepared to accept $b$, so there are no choices for these slots. The remaining $r - s + 1$ slots to eventually receive $c$ must all be distributed, with exactly 1 slot left to the final action of $\langle c \rangle$. This results in $\binom{r-s+1}{1} = r - s + 1$ ideals after the $\langle c^r \rangle$ action. After this choice is made, the remaining choices are predetermined. Hence the coefficient of $\mathcal{M}'$ in $A$ is $r - s + 1$ by Remark 5.1(b). The $m$-tuple $\mathcal{M}'$ was produced in $B$ directly by filling all $r + 1$ slots that accept $c$ first and then filling the remaining $p$ slots that accept $b$. All choices are predetermined, so by Remark 5.1(b) the coefficient of $\mathcal{M}'$ in $B$ is $1 + r - p$ coming from the factor already present. Finally, $\mathcal{M}'$ was also produced in $C$ by choosing the one slot for $b$ among those that do not require $c$ first. There are $p - s$ such slots, and so this results in in $\binom{p-s}{1} = p - s$ ideals after the $\langle b \rangle$ action. Then the remaining choices are predetermined by filling in the $r + 1$ copies of $c$ and the remaining $p - 1$ copies of $b$. So by Remark 5.1(b) the coefficient of $\mathcal{M}'$ in $C$ will be $p - s$. Adding the coefficients for $B$ and $C$, we get $(1 + r - p) + (p - s) = 1 + r - s$, matching the coefficient of $\mathcal{M}'$ in $A$. Hence the base case is proved.

The inductive step is separated into two cases: (i) $1 \leq q < p$ and (ii) $q \geq p$. First assume that $1 \leq q < p$ and that the result holds for all $1 \leq k \leq q$. We indicate when this inductive hypothesis (IH) is used and when Lemma 5.5 is used with references in the right column. Note that $\min\{p, q\} = q$ and that the gluing identity gives $(q + 1)\langle c^{q+1}, b^p, c^r \rangle = \langle c, c^q, b^p, c^r \rangle$, and so we have

$$(q + 1)\langle c^{q+1}, b^p, c^r \rangle = \sum_{k=0}^{q} \binom{q+r-p}{q-k} \langle c, b^{p-k}, c^{q+r}, b^k \rangle \tag{IH}$$

$$= \sum_{k=0}^{q} \binom{q+r-p}{q-k} \left[ (1 + (q+r) - (p-k))\langle b^{p-k}, c^{1+q+r}, b^k \rangle + \langle b^{p-k-1}, c^{1+q+r}, b, b^k \rangle \right] \tag{IH}$$

$$= \sum_{k=0}^{q} \binom{q+r-p}{q-k} (1 + q + r - p + k)\langle b^{p-k}, c^{1+q+r}, b^k \rangle$$

$$+ \sum_{k=0}^{q} \binom{q+r-p}{q-k} (1+k)\langle b^{p-(k+1)}, c^{1+q+r}, b^{k+1} \rangle$$

$$= \sum_{k=0}^{q} \binom{q+r-p}{q-k} (1 + q + r - p + k)\langle b^{p-k}, c^{1+q+r}, b^k \rangle$$

$$+ \sum_{k=1}^{q+1} \binom{q+r-p}{q-k+1} k\langle b^{p-k}, c^{1+q+r}, b^k \rangle$$

$$= \binom{q+r-p}{q}(1 + q + r - p)\langle b^p, c^{(q+1)+r} \rangle$$

$$+ \sum_{k=1}^{q} \left[ \binom{q+r-p}{q-k}(1+q+r-p+k) + \binom{q+r-p}{q-k+1}k \right] \langle b^{p-k}, c^{(q+1)+r}, b^k \rangle$$

$$+ (q+1)\langle b^{p-(q+1)}, c^{(q+1)+r}, b^{q+1} \rangle$$

$$= (q+1)\binom{(q+1)+r-p}{q+1}\langle b^p, c^{(q+1)+r} \rangle \qquad \text{5.5(a)}$$

$$+ \sum_{k=1}^{q}(q+1)\binom{(q+1)+r-p}{(q+1)-k}\langle b^{p-k}, c^{(q+1)+r}, b^k \rangle \qquad \text{5.5(b)}$$

$$+ (q+1)\langle b^{p-(q+1)}, c^{(q+1)+r}, b^{q+1} \rangle$$

$$= (q+1)\sum_{k=0}^{q+1}\binom{(q+1)+r-p}{(q+1)-k}\langle b^{p-k}, c^{(q+1)+r}, b^k \rangle$$

Dividing by $q+1$ completes case (i) where $1 \leq q < p$. Now suppose $q \geq p$ and that the result holds for all $1 \leq k \leq q$. Note that $\min\{p,q\} = p$ and that the gluing identity again gives $(q+1)\langle c^{q+1}, b^p, c^r \rangle = \langle c, c^q, b^p, c^r \rangle$, so we have

$$(q+1)\langle c^{q+1}, b^p, c^r \rangle = \sum_{k=0}^{p}\binom{q+r-p}{q-k}\langle c, b^{p-k}, c^{q+r}, b^k \rangle \qquad \text{(IH)}$$

$$= \sum_{k=0}^{p-1}\binom{q+r-p}{q-k}\langle c, b^{p-k}, c^{q+r}, b^k \rangle + \binom{q+r-p}{q-p}\langle c, c^{q+r}, b^p \rangle$$

$$= \sum_{k=0}^{p-1}\binom{q+r-p}{q-k}\langle c, b^{p-k}, c^{q+r}, b^k \rangle + \binom{q+r-p}{q-p}((q+1)+r)\langle c^{(q+1)+r}, b^p \rangle$$

$$= \sum_{k=0}^{p-1}\binom{q+r-p}{q-k}\left[(1+(q+r)-(p-k))\langle b^{p-k}, c^{((q+1)+r)}, b^k \rangle + \langle b^{p-k-1}, c^{((q+1)+r)}, b, b^k \rangle\right] \qquad \text{(IH)}$$

$$+ \binom{q+r-p}{q-p}((q+1)+r)\langle c^{(q+1)+r}, b^p \rangle$$

$$= \sum_{k=0}^{p}\binom{q+r-p}{q-k}((q+1)+r-p+k)\langle b^{p-k}, c^{(q+1)+r}, b^k \rangle$$

$$+ \sum_{k=1}^{p}\binom{q+r-p}{q-k+1}k\langle b^{p-k}, c^{(q+1)+r}, b^k \rangle$$

$$= \binom{q+r-p}{q}((q+1)+r-p)\langle b^p, c^{(q+1)+r} \rangle$$

$$+ \sum_{k=1}^{p}\left[\binom{q+r-p}{q-k}((q+1)+r-p+k) + \binom{q+r-p}{q-k+1}k\right]\langle b^{p-k}, c^{(q+1)+r}, b^k \rangle$$

$$= (q+1)\binom{(q+1)+r-p}{q+1}\langle b^p, c^{(q+1)+r} \rangle \qquad \text{5.5(a)}$$

$$+ \sum_{k=1}^{p}(q+1)\binom{(q+1)+r-p}{(q+1)-k}\langle b^{p-k}, c^{(q+1)+r}, b^k \rangle \qquad \text{5.5(b)}$$

$$= (q+1) \sum_{k=0}^{p} \binom{(q+1)+r-p}{(q+1)-k} \langle b^{p-k}, c^{(q+1)+r}, b^k \rangle$$

Dividing by $q+1$ completes case (ii) where $q \geq p$. □

## 6 The stackwise spanning theorem

In the next three sections, for notational reasons we ignore operators $\langle c^q \rangle$ if $q = 0$; e.g. $\langle a^2, b^0, c \rangle = \langle a^2, c \rangle$. We likewise suppress poset elements in our standard notation that do not have positive values in a $P$-partition bounded by $m$. For example, if we have the linear extension $v < w < x < y < z$ in Figure 3.1 with $P$-partition values $\psi(v) = 0$, $\psi(w) = 2$, $\psi(x) = 0$, $\psi(y) = 1$, and $\psi(z) = 0$ bounded by 2, we now suppress $v, x,$ and $z$ and write $x_1 = w$ and $x_2 = y$. So $a_1 = c$ and $a_2 = d$ giving the corresponding 2-stackwise operator $\langle a_2^1, a_1^2 \rangle$.

Here we obtain Theorem 6.1, which shows how certain vectors in $\mathcal{V}^m$ can be written as integer linear sums of $m$-stackwise vectors. Their operators have the form $\langle c^q, a_i^{n_i}, \ldots, a_2^{n_2}, a_1^{n_1} \rangle$ where $a_i^{n_i}, \ldots, a_2^{n_2}, a_1^{n_1}$ is an $m$-stackwise word. We say such operators and their words have *augmented stackwise form*.

The proof of Theorem 6.1 inducts on the total word order of Definition 4.7. It relies on key lemmas making use of Propositions 5.2, 5.4, and 5.6. We state these lemmas after stating Theorem 6.1 but delay their proofs until the next section. We now need the full assumption that $P$ is $\Gamma$-colored $d$-complete.

**Theorem 6.1.** *Let $P$ be a $\Gamma$-colored d-complete poset. Fix $m \geq 1$ and a linear extension of $P$. Suppose that $c^q, a_i^{n_i}, \ldots, a_2^{n_2}, a_1^{n_1}$ has augmented stackwise form for $c \in \Gamma$ and $q \geq 1$ and $i \geq 0$ with $n_t > 0$ whenever $1 \leq t \leq i$. The following list encompasses all possibilities for computing $\langle c^q, a_i^{n_i}, \ldots, a_2^{n_2}, a_1^{n_1} \rangle.\{\emptyset, \ldots, \emptyset\}$.*

- *First suppose that $i = 0$.*
  - *(a) If $c$ is not the color of a minimal element of $P$ or if $q > m$, then $\langle c^q \rangle.\{\emptyset, \ldots, \emptyset\} = 0$.*
  - *(b) If $c$ is the color of a minimal element $z$ of $P$ and $q \leq m$, then $\langle c^q \rangle.\{\emptyset, \ldots, \emptyset\} = \binom{m}{q}\{z, \ldots, z, \emptyset, \ldots, \emptyset\}$, where there are $q$ copies of $z$ and $m - q$ copies of $\emptyset$. This is an $m$-stackwise vector.*
- *Now suppose that $i \geq 1$.*
  - *(c) If $c = a_i$ and $q + n_i > m$, then $\langle c^q, a_i^{n_i}, \ldots, a_2^{n_2}, a_1^{n_1} \rangle.\{\emptyset, \ldots, \emptyset\} = 0$.*
  - *(d) If $c^q, a_i^{n_i}, \ldots, a_2^{n_2}, a_1^{n_1}$ is an $m$-stackwise word, then*
    - *(i) If $c = a_i$ then $\langle c^q, a_i^{n_i}, \ldots, a_2^{n_2}, a_1^{n_1} \rangle.\{\emptyset, \ldots, \emptyset\} = \binom{q+n_i}{q}\langle a_i^{q+n_i}, \cdots, a_2^{n_2}, a_1^{n_1} \rangle.\{\emptyset, \ldots, \emptyset\}$, or*
    - *(ii) If $c \neq a_i$ then $\langle c^q, a_i^{n_i}, \ldots, a_2^{n_2}, a_1^{n_1} \rangle.\{\emptyset, \ldots, \emptyset\}$ is an $m$-stackwise vector.*
    
    *In both (i) and (ii) the result is an integer multiple of an $m$-stackwise vector.*
  - *(e) Otherwise $\langle c^q, a_i^{n_i}, \ldots, a_2^{n_2}, a_1^{n_1} \rangle.\{\emptyset, \ldots, \emptyset\}$ is a (possibly zero) integer linear sum of $m$-stackwise vectors of the form $\langle d_h^{o_h}, \ldots, d_2^{o_2}, d_1^{o_1} \rangle.\{\emptyset, \ldots, \emptyset\}$, where $d_h^{o_h}, \ldots, d_2^{o_2}, d_1^{o_1}$ has the same length as and is less than $c^q, a_i^{n_i}, \ldots, a_2^{n_2}, a_1^{n_1}$ in the total word order of Definition 4.7.*

*In all cases the result either vanishes or is an integer linear sum of $m$-stackwise vectors whose corresponding words have the same length as and are less than or equal to $c^q, a_i^{n_i}, \ldots, a_1^{n_1}$ in the total word order.*

The four key lemmas we will need are now given. In each statement the setting of Theorem 6.1 applies; each lemma is used exclusively within the inductive step of the proof of Theorem 6.1 so their proofs use its inductive hypothesis. Note that $I$ in the statements of Lemmas 6.3 and 6.5 are ideals since they correspond to all elements of $P$ whose integer assignments are positive in a $P$-partition bounded by $m$.

**Lemma 6.2.** *Let $c^q, a_i^{n_i}, \ldots, a_1^{n_1}$ be a word of augmented stackwise form and let $e_g^{p_g}, \ldots, e_1^{p_1}$ be an m-stackwise word. Let $f_z^{q_z}, \ldots, f_1^{q_1}$ be such that $c^q, a_i^{n_i}, \ldots, a_1^{n_1}$ and $f_z^{q_z}, \ldots, f_1^{q_1}, e_g^{p_g}, \ldots, e_1^{p_1}$ have the same length. Suppose that there exists some $1 \le h \le p_g + \cdots + p_1$ such that the right initial segments satisfy $R_h(e_g^{p_g}, \ldots, e_1^{p_1}) < R_h(c^q, a_i^{n_i}, \ldots, a_1^{n_1})$. Then $\langle f_z^{q_z}, \ldots, f_1^{q_1}, e_g^{p_g}, \ldots, e_1^{p_1} \rangle.\{\emptyset, \ldots, \emptyset\}$ is a (possibly zero) integer linear sum of m-stackwise vectors of the form $\langle d_h^{o_h}, \ldots, d_2^{o_2}, d_1^{o_1} \rangle.\{\emptyset, \ldots, \emptyset\}$, where $d_h^{o_h}, \ldots, d_2^{o_2}, d_1^{o_1}$ has the same length as and is less than $c^q, a_i^{n_i}, \ldots, a_1^{n_1}$ in the total word order.*

**Lemma 6.3.** *Let $a_i^{n_i}, \ldots, a_2^{n_2}, a_1^{n_1}$ be an m-stackwise word with $i \ge 1$ and $n_t > 0$ for $1 \le t \le i$. Set $I = \{x_i, \ldots, x_2, x_1\}$. Let $c \in \Gamma$ be such that $c \ne a_i$. Suppose that no minimal element of the filter $F = P - I$ has color $c$ and that there are no neighbors in $I$ with colors $c$ and $a_i$. Let $q \ge 1$. Then $\langle c^q, a_i^{n_i}, \ldots, a_1^{n_1} \rangle.\{\emptyset, \ldots, \emptyset\}$ is a (possibly zero) integer linear sum of m-stackwise vectors of the form $\langle d_h^{o_h}, \ldots, d_2^{o_2}, d_1^{o_1} \rangle.\{\emptyset, \ldots, \emptyset\}$, where $d_h^{o_h}, \ldots, d_2^{o_2}, d_1^{o_1}$ has the same length as and is less than $c^q, a_i^{n_i}, \ldots, a_2^{n_2}, a_1^{n_1}$ in the total word order.*

**Lemma 6.4.** *Let $a_i^{n_i}, \ldots, a_2^{n_2}, a_1^{n_1}$ be an m-stackwise word with $i \ge 1$ and $n_t > 0$ for $1 \le t \le i$. Suppose that $x$ with color $c$ is greater than $x_i$ in the linear extension of $P$ and that there is some $x_j$ covered by $x$ such that $q > n_j \ge 1$. Then $\langle c^q, a_i^{n_i}, \ldots, a_1^{n_1} \rangle.\{\emptyset, \ldots, \emptyset\}$ is a (possibly zero) integer linear sum of m-stackwise vectors of the form $\langle d_h^{o_h}, \ldots, d_2^{o_2}, d_1^{o_1} \rangle.\{\emptyset, \ldots, \emptyset\}$, where $d_h^{o_h}, \ldots, d_2^{o_2}, d_1^{o_1}$ has the same length as and is less than $c^q, a_i^{n_i}, \ldots, a_2^{n_2}, a_1^{n_1}$ in the total word order.*

**Lemma 6.5.** *Let $a_i^{n_i}, \ldots, a_2^{n_2}, a_1^{n_1}$ be an m-stackwise word with $n_t > 0$ for $1 \le t \le i$. Set $I = \{x_i, \ldots, x_2, x_1\}$. Let $c \in \Gamma$ be such that $c \ne a_i$. Suppose that no minimal element of the filter $F = P - I$ has color $c$ and that there exists an element $x_j \in I$ with color $c$ that is covered by $x_i$. Then $\langle c^q, a_i^{n_i}, \ldots, a_1^{n_1} \rangle.\{\emptyset, \ldots, \emptyset\}$ is a (possibly zero) integer linear sum of m-stackwise vectors of the form $\langle d_h^{o_h}, \ldots, d_2^{o_2}, d_1^{o_1} \rangle.\{\emptyset, \ldots, \emptyset\}$, where $d_h^{o_h}, \ldots, d_2^{o_2}, d_1^{o_1}$ has the same length as and is less than $c^q, a_i^{n_i}, \ldots, a_2^{n_2}, a_1^{n_1}$ in the total word order.*

*Proof of Theorem 6.1.* We begin by showing that the theorem is true when $i = 0$. If $c$ is not the color of a minimal element of $P$, then $X_c$ vanishes at $\emptyset$ and so $\langle c^q \rangle.\{\emptyset, \ldots, \emptyset\} = 0$. Now suppose $q > m$. The previous case applies if $c$ is still not the color of a minimal element of $P$, so assume $c$ is the color of a minimal element $z$ of $P$. By EC and NA there are no other minimal elements of $P$ of color $c$ and there are no elements covering $z$ of color $c$. So we get $\langle c^q \rangle.\{\emptyset, \ldots, \emptyset\} = \frac{m!(q-m)!}{q!}\langle c^{q-m} \rangle.\{z, \ldots, z\} = 0$ since the next raising operator $X_c$ will vanish at every copy of $z$. Now suppose that $c$ is the color of a minimal element $z$ of $P$ and that $q \le m$. Once again we know there are no other minimal elements of $P$ of color $c$ and no elements covering $z$ of color $c$, so we get

$$\langle c^q \rangle.\{\emptyset, \ldots, \emptyset\} = \frac{1}{q!}\langle c, \ldots, c \rangle.\{\emptyset, \ldots, \emptyset\} = \frac{1}{q!} \cdot \frac{m!}{(m-q)!}\{z, \ldots, z, \emptyset, \ldots, \emptyset\} = \binom{m}{q}\{z, \ldots, z, \emptyset, \ldots, \emptyset\}.$$

(This is a simplified version of Remark 5.1(c).) So if $i = 0$ then either (a) or (b) apply.

We induct on the position of $c^q, a_i^{n_i}, \ldots, a_1^{n_1}$ in the total word order of Definition 4.7. Expand the linear extension of $P$ to $P \cup \Gamma$ where every element of $P$ is greater than every element of $\Gamma$. The lowest word in this order has length 1 with $q = 1$ and $i = 0$. The above argument thus handles the base case.

Now let $c^q, a_i^{n_i}, \ldots, a_1^{n_1}$ have augmented stackwise form and suppose the theorem holds for every augmented stackwise word less than $c^q, a_i^{n_i}, \ldots, a_1^{n_1}$ in the total word order. If $i = 0$, then the above argument shows either (a) or (b) applies, so assume that $i \geq 1$. For notational simplicity, set $b := a_i$ and $p := n_i$ and $x := x_i$. With this new notation the action being considered is $\langle c^q, b^p, a_{i-1}^{n_{i-1}}, \ldots, a_1^{n_1} \rangle.\{\emptyset, \ldots, \emptyset\}$.

First consider the case that $c = b$ and suppose that $q+p > m$. Note that the operator $\langle c^q, c^p, a_{i-1}^{n_{i-1}}, \ldots, a_1^{n_1} \rangle$ concludes with $X_c^{q+p}$. By EC and NA we have $X_c^2.I = 0$ for every ideal $I$ of $P$, so no ideal in any $m$-multiset of ideals arising from this action can accept more than one element of color $c$. Hence $X_c^{q+p}.\{I_1, \ldots, I_m\} = 0$ for every $m$-multiset $\{I_1, \ldots, I_m\}$ of ideals of $P$, and thus (c) holds. Similar reasoning shows the action vanishes if $q > m$ whether or not $c = b$, so for the remainder of the proof we assume $q \leq m$.

Continue to assume $c = b$ and suppose that $q + p \leq m$. If $n_j \geq q + p$ whenever $x_j$ is covered by $x$ in the partial order on $P$ for $1 \leq j \leq i - 1$, then by definition $c^q, c^p, a_{i-1}^{n_{i-1}}, \ldots, a_1^{n_1}$ is an $m$-stackwise word. Then (d)(i) holds once the gluing coefficient is applied. On the other hand, if there exists an $x_j$ covered by $x$ in the partial order on $P$ for some $1 \leq j \leq i - 1$ such that $q + p > n_j$, then $c^q, c^p, a_{i-1}^{n_{i-1}}, \ldots, a_1^{n_1}$ is not an $m$-stackwise word. We apply Lemma 6.4 to the action $\binom{p+q}{p}\langle c^{p+q}, a_{i-1}^{n_{i-1}}, \ldots, a_1^{n_1} \rangle.\{\emptyset, \ldots, \emptyset\}$ to show that (e) holds in this case. (In the statement of Lemma 6.4, replace $q$ with $p + q$ and $i$ with $i - 1$, noting that $i - 1 \geq 1$.)

Now suppose that $c \neq b$ and continue considering $\langle c^q, b^p, a_{i-1}^{n_{i-1}}, \ldots, a_1^{n_1} \rangle.\{\emptyset, \ldots, \emptyset\}$. Let $I = \{x_1, \ldots, x_i\}$ be the ideal of elements whose corresponding $P$-partition integer values are positive and let $F = P - I$. Suppose there is a minimal element $y$ of $F$ with $\kappa(y) = c$.

First suppose that $y < x_i$ in the extension. Then $c^q, b^p, a_{i-1}^{n_{i-1}}, \ldots, a_1^{n_1}$ is not an $m$-stackwise word since it is written in the wrong order for its color groups. In this case, insert $y$ to its correct place in the extension $y < x_k < \cdots < x_i$ for $k \leq i$. Since $y \notin I$, we see that $y$ is incomparable in $P$ to every element $x_k, \ldots, x_i$. Hence $c$ is distant to all of $a_k, \ldots, a_i$ by EC and AC. By Proposition 5.2 it follows that $\langle c^q, a_i^{n_i}, \ldots, a_1^{n_1} \rangle = \langle a_i^{n_i}, \ldots, a_k^{n_k}, c^q, a_{k-1}^{n_{k-1}}, \ldots, a_1^{n_1} \rangle$.

Note that $\{x_1, \ldots, x_{k-1}\}$ is an ideal of $P$. Moreover, since $y$ is minimal in $F$ and is incomparable in $P$ to $x_k, \ldots, x_i$, we see $\{x_1, \ldots, x_{k-1}, y\}$ is an ideal of $P$. Let $h = 1$ if $k = 1$ or $h = n_1 + \cdots + n_{k-1} + 1$ if $k > 1$. Since $\{x_1, \ldots, x_{k-1}, y\}$ is an ideal of $P$, the right initial segment $R_h(c^q, a_{k-1}^{n_{k-1}}, \ldots, x_1^{n_1}) = c, a_{k-1}^{n_{k-1}}, \ldots, a_1^{n_1}$ is an $m$-stackwise word. Both it and $R_h(c^q, a_i^{n_i}, \ldots, a_1^{n_1})$ grow well for $m$ by Lemma 4.3 since they are initial segments of $m$-stackwise words. Then in the notation of Definition 4.7 we have $[R_h(c^q, a_{k-1}^{n_{k-1}}, \ldots, a_1^{n_1})]_\ell = [R_h(c^q, a_i^{n_i}, \ldots, a_1^{n_1})]_\ell$ whenever $\ell < h$ and $[R_h(c^q, a_{k-1}^{n_{k-1}}, \ldots, a_1^{n_1})]_h = y < x_k = [R_h(c^q, a_i^{n_i}, \ldots, a_1^{n_1})]_h$.

Since $c^q, a_{k-1}^{n_{k-1}}, \ldots, a_1^{n_1}$ has augmented stackwise form and is less than the original word $c^q, a_i^{n_i}, \ldots, a_1^{n_1}$ by length, the inductive hypothesis applies and so we can expand $\langle c^q, a_{k-1}^{n_{k-1}}, \ldots, a_1^{n_1} \rangle.\{\emptyset, \ldots, \emptyset\}$ as an integral sum of $m$-stackwise vectors whose words are all less than or equal to $c^q, a_{k-1}^{n_{k-1}}, \ldots, a_1^{n_1}$ in the total word order. Let $\mathbf{d}$ be an $m$-stackwise word appearing in this expansion. By the previous paragraph, we have $R_h(\mathbf{d}) \leq R_h(c^q, a_{k-1}^{n_{k-1}}, \ldots, a_1^{n_1}) < R_h(c^q, a_i^{n_i}, \ldots, a_1^{n_1})$. We then apply Lemma 6.2 to the expansion of $\langle a_i^{n_i}, \ldots, a_k^{n_k}, \mathbf{d} \rangle.\{\emptyset, \ldots, \emptyset\}$ to show (e) is satisfied in this case.

Continue to suppose that $c \neq b$ and that there is a minimal element $y$ of $F$ with $\kappa(y) = c$, but now suppose that $y > x_i$ in the extension. If $y$ is minimal in $P$, then $c^q, b^p, a_{i-1}^{n_{i-1}}, \ldots, a_1^{n_1}$ is an $m$-stackwise word by definition and d(ii) applies, so suppose $y$ is not minimal in $P$. Since $y$ is minimal in $F$, we note it only covers elements of $I$. If $q \leq \min\{n_j \mid y \text{ covers } x_j\}$, then $c^q, b^p, a_{i-1}^{n_{i-1}}, \ldots, a_1^{n_1}$ is an $m$-stackwise word by definition and d(ii) applies. If $q > \min\{n_j \mid y \text{ covers } x_j\}$, then Lemma 6.4 can be used again to show that (e) holds. (In the statement of Lemma 6.4, replace $x$ with $y$.)

Continue to assume that $c \neq b$ and now suppose that no such minimal element $y$ of $F$ with color $c$ exists. Then $c^q, b^p, a_{i-1}^{n_{i-1}}, \ldots, a_1^{n_1}$ is not an $m$-stackwise word. If there are no neighbors in $I$ with colors $b$ and $c$, then Lemma 6.3 shows that (e) holds in this case. Now suppose that $b$ and $c$ are colors of neighbors in $I$. By NA we know that $b$ and $c$ are adjacent. We claim that an element of color $c$ must be covered by $x = x_i$. Note that $P_c \cap I \neq \emptyset$ by assumption, so let $u$ be the maximal element of $P_c \cap I$. By AC and since $u \in I$ we see that $u < x$ in $P$. If $u$ is maximal in $P_c$, then $u \to x$ by NA and UCB1. If $u$ is not maximal in $P_c$, then let $v$ be minimal in $P_c \cap F$. By EC we know that $u < v$ are consecutive elements of color $c$ and by AC we know that $x \in [u, v]$. Either $u \to x$ or $x \to v$ by NA and ICE2. Suppose that $x \to v$. Then since $x$ is maximal in $I$, we must have $v \in F$. Note our assumption that $v$ is not minimal in $F$, so there is some element $w \to v$ with $w \in F$. We know $w \neq x$ since $x \in I$. By NA and AC we must have $w \in [u, v]$ as well. But then by NA and ICE2 we must have $u \to w$ and $u \to x$. In either case we have $u \to x$, so the claim is proved. Thus Lemma 6.5 shows that (e) holds in this case.

This represents all scenarios, and so the proof is complete once the four lemmas have been proved. □

## 7 Proofs of the four key stackwise spanning lemmas

To complete the proof of Theorem 6.1, we prove here Lemmas 6.2, 6.3, 6.4, and 6.5. These lemmas are all situated within the inductive step of the proof of Theorem 6.1. We state here the stackwise spanning theorem inductive hypothesis (SSTIH) for reference, which we invoke in all four proofs.

**SSTIH**: Let $c^q, a_i^{n_i}, \ldots, a_1^{n_1}$ be a word of augmented stackwise form as in the statement of Theorem 6.1. Assume that this theorem holds for every augmented stackwise word less than $c^q, a_i^{n_i}, \ldots, a_1^{n_1}$ in the total word order of Definition 4.7.

Before their proofs, we remark that if **a** and **b** are words of the same length $k$ and the initial portions of length $h$ satisfy $R_h(\mathbf{a}) < R_h(\mathbf{b})$ for some $1 \leq h \leq k$, then $\mathbf{a} < \mathbf{b}$ in the lexicographic total word order. See Section 6 for the statements of these four lemmas.

*Proof of Lemma 6.2.* By word length (if $z > 1$) or the remark above (if $z = 1$) we see that $f_1^{q_1}, e_g^{p_g}, \ldots, e_1^{p_1}$ is less than $c^q, a_i^{n_i}, \ldots, a_1^{n_1}$ in the total word order. This word $f_1^{q_1}, e_g^{p_g}, \ldots, e_1^{p_1}$ has augmented stackwise form, so the SSTIH applies. Hence one of (a)–(e) holds and we may expand $\langle f_1^{q_1}, e_g^{p_g}, \ldots, e_1^{p_1} \rangle.\{\emptyset, \ldots, \emptyset\}$ as a (possibly zero) integer linear sum of $m$-stackwise vectors, each of whose corresponding words are less than or equal to $f_1^{q_1}, e_g^{p_g}, \ldots, e_1^{p_1}$ in the total word order but of the same length. Suppose **d** is one of the $m$-stackwise words from this expansion so that $\mathbf{d} \leq f_1^{q_1}, e_g^{p_g}, \ldots, e_1^{p_1}$. Then $f_2^{q_2}, \mathbf{d} \leq f_2^{q_2}, f_1^{q_1}, e_g^{p_g}, \ldots, e_1^{p_1} < c^q, a_i^{n_i}, \ldots, a_1^{n_1}$, and $f_2^{q_2}, \mathbf{d}$ has augmented stackwise form and so the SSTIH applies. This reasoning can be repeated until we have expanded $\langle f_z^{q_z}, \ldots, f_1^{q_1}, e_g^{p_g}, \ldots, e_1^{p_1}\rangle.\{\emptyset, \ldots, \emptyset\}$. □

*Proof of Lemma 6.3.* First suppose that $P_c \cap I = \emptyset$. By Lemma 4.3, the ideal $I$ is the first ideal in the $m$-flag corresponding to the $m$-stackwise word $a_i^{n_i}, \ldots, a_1^{n_1}$. Hence by Remark 4.5 only elements from $I$ can be added to any ideal in any $m$-multiset of ideals in the expansion of $\langle a_i^{n_i}, \ldots, a_1^{n_1}\rangle.\{\emptyset, \ldots, \emptyset\}$. Since $F$ has no minimal elements of color $c$, we conclude that $\langle c^q, a_i^{n_i}, \ldots, a_1^{n_1}\rangle.\{\emptyset, \ldots, \emptyset\} = 0$.

Now suppose that $P_c \cap I \neq \emptyset$. Set $y = x_i$ and $b = a_i$ and $p = n_i$ so that the operator we are considering is $\langle c^q, a_i^{n_i}, \ldots, a_1^{n_1}\rangle = \langle c^q, b^p, a_{i-1}^{n_{i-1}}, \ldots, a_1^{i_1}\rangle$. For reasons similar to those given in the above paragraph, since there are no minimal elements of $F$ of color $c$ only elements in $I$ can appear in any ideal of any $m$-multiset of ideals in the expansion of $\langle c^q, b^p, a_{i-1}^{n_{i-1}}, \ldots, a_1^{n_1}\rangle.\{\emptyset, \ldots, \emptyset\}$. Since no elements of colors $b$ and $c$ are neighbors in $I$, we apply Proposition 5.2 to the poset $I$ to get

$$\langle c^q, b^p, a_{i-1}^{n_{i-1}}, \ldots, a_1^{n_1}\rangle.\{\emptyset, \ldots, \emptyset\} = \langle b^p, c^q, a_{i-1}^{n_{i-1}}, \ldots, a_1^{n_1}\rangle.\{\emptyset, \ldots, \emptyset\}$$

which holds when viewed in $P$ as well. We are done if this action vanishes, so assume it does not vanish. Since $c^q, a_{i-1}^{n_{i-1}}, \ldots, a_1^{n_1}$ is shorter and hence less than $c^q, b^p, a_{i-1}^{n_{i-1}}, \ldots, a_1^{n_1}$ in the total word order and has augmented stackwise form, we apply the SSTIH to expand $\langle c^q, a_{i-1}^{n_{i-1}}, \ldots, a_1^{n_1}\rangle.\{\emptyset, \ldots, \emptyset\}$ as an integral linear sum of $m$-stackwise vectors whose corresponding words are less than or equal to $c^q, a_{i-1}^{n_{i-1}}, \ldots, a_1^{n_1}$ in the total word order and have the same length as $c^q, a_{i-1}^{n_{i-1}}, \ldots, a_1^{n_1}$. Let $\mathbf{d}$ be an $m$-stackwise word in this expansion. Even though $b^p, \mathbf{d}$ has augmented stackwise form, we must be careful to check its position relative to $c^q, b^p, a_{i-1}^{n_{i-1}}, \ldots, a_1^{n_1}$ in the total word order. We do this and apply Lemma 6.2.

Again using Lemma 4.3 and Remark 4.5, every ideal in every $m$-multiset of ideals in the expansion of $\langle a_{i-1}^{n_{i-1}}, \ldots, a_1^{n_1}\rangle.\{\emptyset, \ldots, \emptyset\}$ is contained in the ideal $I' = I - \{y\}$. Let $z$ be an element of color $c$ added to an ideal during the action $\langle c^q \rangle$ on the right hand side operator above. Suppose for a contradiction that $z$ is greater than $y$ in the extension. Then $z \notin I'$ and so $z$ is minimal in the filter $F' = P - I'$. Since $y$ is also minimal in $F'$, we see $y$ and $z$ are incomparable in $P$, contradicting AC. Hence $z < y$ in the extension.

Set $k' = n_{i-1} + \cdots + n_1 + 1$. We claim the initial segment relation

$$R_{k'}(c^q, a_{i-1}^{n_{i-1}}, \ldots, a_1^{n_1}) < R_{k'}(c^q, b^p, a_{i-1}^{n_{i-1}}, \ldots, a_1^{n_1}).$$

The word $b, a_{i-1}^{n_{i-1}}, \ldots, a_1^{n_1}$ grows well for $m$ by Lemma 4.3 since it is the initial segment of an $m$-stackwise word, so we have $[b, a_{i-1}^{n_{i-1}}, \ldots, a_1^{n_1}] = y, x_{i-1}^{n_{i-1}}, \ldots, x_1^{n_1}$ in the notation of Definition 4.7. We observe that $c, a_{i-1}^{n_{i-1}}, \ldots, a_1^{n_1}$ may or may not grow well for $m$ but note that the initial right segment of length $k' - 1$ does grow well for $m$ since it is an $m$-stackwise word. So either $[c, a_{i-1}^{n_{i-1}}, \ldots, a_1^{n_1}] = c, x_{i-1}^{n_{i-1}}, \ldots, x_1^{n_1}$ if this word does not grow well for $m$ or $[c, a_{i-1}^{n_{i-1}}, \ldots, a_1^{n_1}] = z, x_{i-1}^{n_{i-1}}, \ldots, x_1^{n_1}$ for some $z$ as in the paragraph above if it does. Since both $c$ and $z$ are less than $y$ in the total order on $P \cup \Gamma$, this proves the claim.

We have now shown that $R_{k'}(\mathbf{d}) \leq R_{k'}(c^q, a_{i-1}^{n_{i-1}}, \ldots, a_1^{n_1}) < R_{k'}(c^q, b^p, a_{i-1}^{n_{i-1}}, \ldots, a_1^{n_1})$. Using $h = k'$ we apply Lemma 6.2 to the expansion of $\langle b^p, \mathbf{d}\rangle.\{\emptyset, \ldots, \emptyset\}$ to show (e) is satisfied. □

*Proof of Lemma 6.4.* Set $y = x_j$ and $b = a_j$ and $p = n_j$. So we have $q > p \geq 1$ and the action we are considering is $\langle c^q, a_i^{n_i}, \ldots, a_1^{n_1} \rangle.\{\emptyset, \ldots, \emptyset\} = \langle c^q, a_i^{n_i}, \ldots, b^p, \ldots, a_1^{n_1} \rangle.\{\emptyset, \ldots, \emptyset\}$.

First suppose that there is some element covered by $x$ which follows $y$ in the extension and is represented in $a_i^{n_i}, \ldots, a_1^{n_1}$. Let $z$ be the greatest element in the extension with this property, and suppose $z$ has color $f$ with multiplicity $r \geq 1$. Note that $y$ and $z$ are incomparable in $P$ because they are both covered by $x$. In this case the operator becomes $\langle c^q, a_i^{n_i}, \ldots, f^r, \ldots, b^p, \ldots, a_1^{n_1} \rangle$. There is no guarantee that $q > r$ and so we cannot replace $y$ with $z$ without loss of generality.

We first claim we can move the $b^p$ action to the left of $f^r$ without affecting the result of the corresponding operator. We consider all groups of colors between $b^p$ and $f^r$ in the word. These correspond to elements that are greater than $y$ but less than $z$ in the extension. Such an element must be (i) greater than $y$ and incomparable to $z$ in $P$, (ii) incomparable to $y$ and less than $z$ in $P$, or (iii) incomparable to both $y$ and $z$ in $P$. It cannot be greater than $y$ and less than $z$ in $P$ since $y$ and $z$ are incomparable in $P$. If an element $u$ in group (i) is less in the extension than an element $v$ in group (ii), then $u$ and $v$ must be incomparable in $P$ or else we would have $u < v < z$ in $P$, violating that $u$ and $z$ are incomparable in $P$. For similar reasons, elements in group (i) that are less than elements in group (iii) in the extension are incomparable in $P$. Incomparable elements have distant colors by EC and AC and groups of distant colors can be moved past each other freely by Proposition 5.2. So color groups corresponding to elements in group (i) can be kept in the same order within group (i) but moved to the left of color groups corresponding to elements in groups (ii) and (iii) and then moved past $f^r$. Then color groups corresponding to elements in groups (ii) and (iii) can be kept in the same order collectively but moved to the right of $b^p$, so that $f^r$ and $b^p$ are consecutive in the operator. Finally, since $y$ and $z$ are incomparable in $P$ and thus have distant colors, we can switch the order of the color groups $f^r$ and $b^p$. Thus we have used Proposition 5.2 to translate the operator of the form $\langle c^q, a_i^{n_i}, \ldots, f^r, (i),(ii),(iii), b^p, \ldots, a_1^{n_1} \rangle$ into the equivalent operator of the form $\langle c^q, a_i^{n_i}, \ldots, (i), b^p, f^r, (ii), (iii), \ldots, a_1^{n_1} \rangle$ where (i), (ii), and (iii) represent all color groups corresponding to elements in groups (i), (ii), and (iii) respectively, kept in their original orders except for the noted switches.

Consider all color groups between $c^q$ and $b^p$ in the newly arranged operator. They correspond to elements that are either in group (i) or between $z$ and $x$ in the extension and represented in $a_i^{n_i}, \ldots, a_1^{n_1}$. Consider such an element $w$. If $w$ is in group (i) above, then $y < w$ in $P$. We also have $w < x$ in the linear extension, so $w$ must be incomparable to $x$ in $P$ since $y \to x$ in $P$. Hence $w$ has color distant to $c$. Now suppose $z < w < x$ in the extension for some $w$ represented in $a_i^{n_i}, \ldots, a_1^{n_1}$. If $w$ is not less than $x$ in $P$, then it is incomparable to $x$ and hence has color distant to $c$. Continue to assume $z < w < x$ in the extension and now suppose that $w < x$ in $P$. If $w$ has color $c$, then since $z \to x$ we have $w < z$ in $P$ by NA and AC, contradicting that $z < w$ in the extension. Hence $w$ does not have color $c$. Now suppose $w$ has color adjacent to $c$. Since $w < x$ in $P$ there is some element $u \to x$ with $w < u$ in $P$. Since $z$ is the maximal element in the extension covered by $x$ in $P$ represented in $a_i^{n_i}, \ldots, a_1^{n_1}$, it must be that $u$ is not represented in $a_i^{n_i}, \ldots, a_1^{n_1}$. Hence $u \neq z$ and $u \neq w$. By NA we know that $u$ and $z$ have colors adjacent to $c$. If there is some element $v$ less than $x$ in $P$ with color $c$, then by NA, AC, and I2A we would have $v \to u$ and $v \to z$ and $w$ is not in the interval $[v, x]$. Then by AC we would have $w < v < z$ in $P$, contradicting that $z < w$ in the extension. Hence no such element $v$ exists. Let $I = \{x_i, \ldots, x_1\}$ be the ideal corresponding to the $P$-partition bounded by $m$ giving $a_i^{n_i}, \ldots, a_1^{n_1}$. We have just shown that $x$ is minimal in $P_c$ and hence $I$ contains no elements of color $c$. Since $u \to x$ and $u \notin I$, we see that $F = P - I$ does not contain a minimal element of color $c$. Thus we can apply Lemma 6.3 to finish the proof after first putting the operator back in its original order. So in all cases either Lemma 6.3 finishes the proof or $w$ has color distant to $c$. We will assume the latter for all such elements $w$. Hence every color appearing in a color group between $c^q$ and $b^p$ in the new form of the operator is distant to $c$, so by Proposition 5.2 we can move $c^q$ to the right until it is next to $b^p$. Hence our original augmentation operator can be converted again to an equivalent form $\langle a_i^{n_i}, \ldots, (i), c^q, b^p, f^r, (ii), (iii), a_{j-1}^{n_{j-1}}, \ldots, a_1^{n_1} \rangle$.

Since $q > p \geq 0$ we now apply Proposition 5.4 to the operator $\langle c^q, b^p \rangle$ sitting in the middle of the above operator to get the final form of our operator identity

$$\langle c^q, a_i^{n_i}, \ldots, f^r, \ldots, b^p, \ldots, a_1^{n_1} \rangle = \sum_{k=0}^{q-1} (-1)^k \langle a_i^{n_i}, \ldots, (i), c^{q-1-k}, b^p, c^{k+1}, f^r, (ii), (iii), a_{j-1}^{n_{j-1}}, \ldots, a_1^{n_1} \rangle$$

The initial right segment $\mathbf{f} := f^r, (ii), (iii), a_{j-1}^{n_{j-1}}, \ldots, a_1^{n_1}$ is the same in each of the resulting operators on the right hand side of the above identity. Note that this is the same initial segment as $c^q, a_i^{n_i}, \ldots, a_1^{n_1} = c^q, a_i^{n_i}, \ldots, f^r, (i), (ii), (iii), b^p, a_{j-1}^{n_{j-1}}, \ldots, a_1^{n_1}$ except for shifting the color groups from elements in group (i) as well as $b^p$ to the left of $f^r$. Now let $I$ be the ideal of $P$ consisting of elements corresponding to the initial segment $f^r, (i), (ii), (iii), b^p, a_{j-1}^{n_{j-1}}, \ldots, a_1^{n_1}$ of the original $m$-stackwise word and let $J$ be the set of elements corresponding to the colors in $\mathbf{f}$. Note that $J \subseteq I$ and that $I - J$ consists of $y$ and every element in group (i). If $u \in I$ with $u > y$ in $P$, then $u$ is incomparable to $z$ and hence $u$ is in (i). So $I - J$ simply consists of $y$ and every element above $y$ in $I$ and hence $J$ is an ideal. Since the multiplicities leading to $a_i^{n_i}, \ldots, a_1^{n_1}$ come from a $P$-partition bounded by $m$ and $\mathbf{f}$ has the same multiplicities for the remaining poset elements, it also corresponds to a $P$-partition bounded by $m$. Hence $\mathbf{f}$ is an $m$-stackwise word.

We now consider $\langle c^g, \mathbf{f} \rangle . \{\emptyset, \ldots, \emptyset\}$ for $1 \leq g \leq q$ from the right hand side above. First assume that $j = 1$ so that $y$ is a minimal element of $P$ of color $b$. Colors appearing in $\mathbf{f} = f^r, (ii), (iii)$ are all distant to $b$ by EC and AC since they come from groups (ii) or (iii) or are $f$. These colors do not include $c$ since $b$ and $c$ are adjacent by NA. So the color group $c^g$ is the first appearance of either $b$ or $c$ in $c^g, \mathbf{f}$. No ideal in any $m$-multiset of ideals in the expansion of $\langle \mathbf{f} \rangle . \{\emptyset, \ldots, \emptyset\}$ contains the minimal element $y$ and so no ideal can accept $x$, the minimal element of color $c$ in $P$. Hence $\langle c^g, \mathbf{f} \rangle . \{\emptyset, \ldots, \emptyset\} = 0$. Now suppose that $j \geq 2$.

The word $\mathbf{f}$ is an $m$-stackwise word and so $c^g, \mathbf{f}$ has augmented stackwise form. It has length less than $q + n_i + \cdots + n_1$ since at least $p$ from $b^p$ is missing from $c^g, \mathbf{f}$. So $c^g, \mathbf{f} < c^q, a_i^{n_i}, \ldots, a_1^{n_1}$ in the total word order and hence the SSTIH applies. We note that $y$ is a minimal element of color $b$ in $P - J$. By AC we know that $y$ must be comparable to all elements of color $c$ in $P$, so no element minimal in $P - J$ can have color $c$. Thus $c^g, \mathbf{f}$ is not an $m$-stackwise word. Hence (a) or (c) or (e) from Theorem 6.1 applies for the expansion of $\langle c^g, \mathbf{f} \rangle . \{\emptyset, \ldots, \emptyset\}$. Both (a) and (c) result in a vanishing action, so for the remainder of this proof we will assume that (e) applies.

Let $\mathbf{d}$ be an $m$-stackwise word from the expansion of $\langle c^g, \mathbf{f} \rangle . \{\emptyset, \ldots, \emptyset\}$ from (e) so that $\mathbf{d} < c^g, \mathbf{f}$ in the total word order. Let $k' = n_{j-1} + \cdots + n_1$. Suppose first that $R_{k'}(\mathbf{d}) \neq R_{k'}(c^g, \mathbf{f})$. Since $\mathbf{d} < c^g, \mathbf{f}$, we must have $R_{k'}(\mathbf{d}) < R_{k'}(c^g, \mathbf{f})$. Furthermore, $R_{k'}(c^g, \mathbf{f}) = R_{k'}(c^q, a_i^{n_i}, \ldots, a_1^{n_1})$. Hence $\mathbf{d}$ is an $m$-stackwise word with $R_{k'}(\mathbf{d}) < R_{k'}(c^q, a_i^{n_i}, \ldots, a_1^{n_1})$. Using $h = k'$ we apply Lemma 6.2 to $\langle a_i^{n_i}, \ldots, (i), c^{q-g}, b^p, \mathbf{d} \rangle.\{\emptyset, \ldots, \emptyset\}$ to show (e) is satisfied for this kind of word $\mathbf{d}$.

Now suppose that $R_{k'}(\mathbf{d}) = R_{k'}(c^g, \mathbf{f}) = a_{j-1}^{n_{j-1}}, \ldots, a_1^{n_1}$. There are $g$ extra copies of $c$ not in $\mathbf{f}$ that must appear somewhere in $\mathbf{d}$ (see Remark 3.2), and we claim they must appear at $a_{j-1}$. They cannot appear in the initial segment $R_{k'}(\mathbf{d})$ since $R_{k'}(\mathbf{d}) = R_{k'}(c^g, \mathbf{f})$. Since all remaining colors to the left of $a_{j-1}$ in $\mathbf{f}$ are distant to $b$ (using the same reasoning as the $j = 1$ case three paragraphs above) and $\mathbf{d}$ and $c^g, \mathbf{f}$ have the same color census, the color $b$ does not appear in $\mathbf{d}$ past the initial segment $R_{k'}(\mathbf{d})$. Thus the element $y$ is not in the ideal $K$ consisting of elements corresponding to $\mathbf{d}$. Since $b$ and $c$ are adjacent, we note by AC that $y$ would have to be in $K$ if any element of color $c$ greater than $y$ in the extension is in $K$. Hence these extra $g$ copies of the color $c$ cannot appear to the left of $a_{j-1}$ in $\mathbf{d}$. Thus they must appear at $a_{j-1}$, as claimed, and add to the total multiplicity at $a_{j-1}$. Set $k'' = k' + 1$ and note that $R_{k''}(\mathbf{d}) = a_{j-1}^{n_{j-1}+1}, a_{j-2}^{n_{j-2}}, \ldots, a_1^{n_1}$ and $R_{k''}(c^q, a_i^{n_i}, \ldots, a_1^{n_1}) = a_j, a_{j-1}^{n_{j-1}}, a_{j-2}^{n_{j-2}}, \ldots, a_1^{n_1}$. Both of these initial segments grow well for $m$ by Lemma 4.3(b) since they are the initial portions of the $m$-stackwise words $\mathbf{d}$ and $a_i^{n_i}, \ldots, a_1^{n_1}$, respectively. So in the total word order of Definition 4.7 we have $R_{k''}(\mathbf{d}) < R_{k''}(c^q, a_i^{n_i}, \ldots, a_1^{n_1})$ since $[R_{k''}(\mathbf{d})]_\ell = [R_{k''}(c^q, a_i^{n_i}, \ldots, a_1^{n_1})]_\ell$ for all $1 \leq \ell \leq k'$ but $[R_{k''}(\mathbf{d})]_{k''} = x_{j-1} < x_j = [R_{k''}(c^q, a_i^{n_i}, \ldots, a_1^{n_1})]_{k''}$. Using $h = k''$ we apply Lemma 6.2 to the expansion of $\langle a_i^{n_i}, \ldots, (i), c^{q-g}, b^p, \mathbf{d} \rangle.\{\emptyset, \ldots, \emptyset\}$ to show (e) is satisfied for this kind of word $\mathbf{d}$. This completes the case where there is some element covered by $x$ which follows $y$ in the extension and is represented in $a_i^{n_i}, \ldots, a_1^{n_1}$.

If there is no such element so that $y$ is the greatest element in the extension covered by $x$ that is represented in $a_i^{n_i}, \ldots, a_1^{n_1}$, then a similar proof proceeds as above without the considerations given to $f^r$ and the color groups corresponding to elements in groups (i), (ii), and (iii) in $P$. □

*Proof of Lemma 6.5.* This proof has a similar structure as the proof of Lemma 6.4. We first show how to rearrange color groups in the operator being considered without changing its output. Let $x = x_i$, $b = a_i$, and $p = n_i$ and let $y = x_j$, $c = a_j$, and $r = n_j$. The operator being considered is $\langle c^q, b^p, a_{i-1}^{n_{i-1}}, \ldots, c^r, \ldots, a_1^{n_1} \rangle$.

First suppose that there is some element covered by $x$ which follows $y$ in the extension. Let $z$ be the greatest element in the extension with this property. Note that $z$ will not have color $c$ so we cannot replace $y$ with $z$ without loss of generality. Note $z$ is represented in $a_i^{n_i}, \ldots, a_1^{n_1}$ since it is covered by $x$ and these multiplicities arise from a $P$-partition bounded by $m$. Suppose $z$ has color $f$ and multiplicity $s$. So the original word is $c^q, b^p, a_{i-1}^{n_{i-1}}, \ldots, f^s, \ldots, c^r, \ldots, a_1^{n_1}$.

Suppose for a contradiction that there is an element $w$ represented in $a_i^{n_i}, \ldots, a_1^{n_1}$ with $y < w < z$ in the extension and with $y < w$ in $P$. Assume that $y \to w$ (or else replace it with such a $w$ that has the same properties). Then $w \neq x$ since $w < z < x$ in the extension. Since $y \to w$ and $y \to x$, by NA and UCB1 there must be an element $v$ of color $c$ greater than $y$. By NA, AC, and ICE2 we must have $w \to v$ and $x \to v$, and $v$ cannot cover any other elements in $P$. Note that both $w$ and $x$ are in $I$. Hence $v$ is minimal in $F = P - I$. This contradicts the assumption that $F$ has no minimal elements of color $c$. Hence no such elements $w$ exist.

Now we consider all color groups between $c^r$ and $f^s$ in the word $c^q, b^p, a_{i-1}^{n_{i-1}}, \ldots, f^s, \ldots, c^r, \ldots, a_1^{n_1}$. These correspond to elements that are greater than $y$ and less than $z$ in the extension. By the previous paragraph, such an element $w$ cannot be greater than $y$ in $P$ so it must be (i) incomparable to $y$ and less than $z$ in $P$, or (ii) incomparable to both $y$ and $z$ in $P$. So the word we are considering has the form $c^q, b^p, a_{i-1}^{n_{i-1}}, \ldots, f^s, (i), (ii), c^r, \ldots, a_1^{n_1}$. Since $y$ is incomparable to elements in groups (i) and (ii) and also to $z$, by EC and AC all of the colors in $f^s$, (i), (ii) are distant to $c$. So by Proposition 5.2 we can move $c^r$ to the left past $f^s$ without affecting the resulting augmentation operator. This obtains the equivalent operator $\langle c^q, b^p, a_{i-1}^{n_{i-1}}, \ldots, c^r, f^s, (i), (ii), a_{j-1}^{n_{j-1}}, \ldots, a_1^{n_1} \rangle$.

Now suppose $z < w < x$ in the extension where $w$ is represented in $a_i^{n_i}, \ldots, a_1^{n_1}$. Then $w$ is incomparable to $x$ in $P$ because $z$ is the maximal element in the extension covered by $x$ in $P$. Hence the colors of every such element $w$ are distant to $b$, and so by Proposition 5.2 the color group $b^p$ can be moved to the right until it is next to the new position of $c^r$. This obtains the equivalent operator $\langle c^q, a_{i-1}^{n_{i-1}}, \ldots, b^p, c^r, f^s, (i), (ii), a_{j-1}^{n_{j-1}}, \ldots, a_1^{n_1} \rangle$.

Suppose for a contradiction that there is some element $u = x_k$ such that $z < u < x$ in the extension with $e = \kappa(u)$ adjacent to $c$. Take $u$ to be the maximal element in the extension with these properties. We claim there is an element $y'$ covered by $u$ with $\kappa(y') = c$. Let $J = \{x_k, \ldots, x_2, x_1\}$ and note $J$ is an ideal of $P$. Note that $P_c \cap J \neq \emptyset$ since $y \in P_c \cap J$. So let $y'$ be the maximal element of $P_c \cap J$. Since $u$ is maximal in $J$, by AC we know $y' < u$ in $P$. If the filter $P - J$ has no elements of color $c$, then $y' \to u$ by NA and UCB1. Otherwise $P - J$ does have an element of color $c$. We can take the minimal such element $y''$ and note that $y' < y''$ are consecutive elements of color $c$ in $P$. By AC we have $u \in [y', y'']$. Either $y' \to u$ or $u \to y''$ by NA, AC, and ICE2. Suppose that $u \to y''$ in $P$. If $y'' \in I$, then $y'' < x$ in $P$ by AC since $c$ and $b$ are adjacent. This would imply that $y'' \leq z$ in the extension by the preceding paragraph since $y''$ has color $c$ (not distant to $b$) and is represented in $a_i^{n_i}, \ldots, a_1^{n_1}$. But we cannot have $y'' \leq z$ in the extension since $z < u$ in the extension and $u \to y''$ in $P$. Hence $y'' \notin I$. Since $y'' \in F$ and $F$ contains no minimal elements of color $c$, there must be some $v \in F$ with $v \to y''$. Note that $v \neq u$ since $v \in F$ and $u \in I$. So $u \to y''$ and $v \to y''$. Hence by NA, AC, and ICE2 we must have $u, v \in [y', y'']$ with $y' \to u$ and $y' \to v$. We still have $y' \to u$, so our claim is proved. By EC we know $y$ and $y'$ are comparable in $P$. If $y' < y$, then since $y' \to u$ we have $u \in [y', y]$ by NA and AC. This would imply $u < y$ in the extension, which is a contradiction since $y < z < u$ in the extension. If $y < y'$, then by AC we would have $x \in [y, y']$. This would imply $x < y'$, meaning that $u < x < y'$ in the extension, which is a contradiction since $y' \to u$ in $P$. Hence we must have $y = y'$. This shows that $y \to x$ and $y \to u$ in $P$. Note that $x \neq u$ since $u < x$ in the extension by assumption. Hence by NA, AC, and UCB1 there is an element $t'$ of color $c$ greater than $x$ and $u$ in $P$. By NA and ICE2 we have $x \to t'$ and $u \to t'$. Moreover, by NA, AC, and ICE2 we know $t'$ cannot cover any other elements in $P$. Since $x < t'$ in $P$, we see $t' \in F$. Thus $t'$ is a minimal element of $F$ of color $c$, which contradicts the original assumption. So we see there is no $u$ with $z < u < x$ in the extension with color adjacent to $c$.

There is also no element $u$ with $z < u < x$ in the extension with $\kappa(u) = c$. This is because $y$ is maximal in $I$ with color $c$ since $y$ is covered by $x$. Hence the color group $c^q$ can be moved to the right by Proposition 5.2 until it is next to $b^p$ to obtain the equivalent operator $\langle a_{i-1}^{n_{i-1}}, \ldots, c^q, b^p, c^r, f^s, (i), (ii), a_{j-1}^{n_{j-1}}, \ldots, a_1^{n_1} \rangle$. Note that $r - p \geq 0$ because $y \to x$ in $P$ and the multiplicities $r$ for $y$ and $p$ for $x$ come from a $P$-partition bounded by $m$. Since this operator is equivalent to the original operator, we use Proposition 5.6 to obtain

$$\langle c^q, b^p, a_{i-1}^{n_{i-1}}, \ldots, c^r, \ldots, a_1^{n_1}\rangle = \sum_{g=0}^{\min\{p,q\}} \binom{q+r-p}{q-g} \langle a_{i-1}^{n_{i-1}}, \ldots, b^{p-g}, c^{q+r}, b^g, f^s, (i), (ii), a_{j-1}^{n_{j-1}}, \ldots, a_1^{n_1}\rangle$$

The rest of the proof will resemble the end of the proof of Lemma 6.4. We will show that every operator in the expansion above can be written as an integral sum of $m$-stackwise vectors by applying Lemma 6.2.

Notice each of the above operators on the right begins with the word $\mathbf{f} = f^s, (i), (ii), a_{j-1}^{n_{j-1}}, \ldots, a_1^{n_1}$. Consider the set $K$ of elements of $P$ corresponding to $\mathbf{f}$. Then $K$ consists of the ideal corresponding to the positive $P$-partition elements up through $z$ in the extension, but with $y$ removed. We showed earlier that all colors in $f^s, (i), (ii)$ are distant to $c$. If $y \to v$ for some $v$ represented in $f^s, (i), (ii)$, then it would have adjacent color by NA. Hence $y$ is maximal in the above mentioned ideal and so $K$ is this ideal with one of its maximal elements removed. Thus $K$ is also an ideal. Since the multiplicities leading to $a_i^{n_i}, \ldots, a_1^{n_1}$ come from a $P$-partition bounded by $m$ and $\mathbf{f}$ has the same multiplicities for the remaining poset elements, it also corresponds to a $P$-partition bounded by $m$. Hence $\mathbf{f}$ is an $m$-stackwise word.

First suppose that $g = 0$ so that the word we are considering has initial segment $c^{q+r}, \mathbf{f}$. In this case move $c^{q+r}$ back to the right past $f^s, (i), (ii)$ to obtain a word that starts with $c^{q+r}, a_{j-1}^{n_{j-1}}, \ldots, a_1^{n_1}$. If $q + r > m$, then the operator of this word on $\{\emptyset, \ldots, \emptyset\}$ vanishes, so suppose $q + r \leq m$. First suppose $c^{q+r}, a_{j-1}^{n_{j-1}}, \ldots, a_1^{n_1}$ is an $m$-stackwise word. In this case, move the color group $b^p$ back to the left so that the full word is $b^p, a_{i-1}^{n_{i-1}}, \ldots, c^{q+r}, a_{j-1}^{n_{j-1}}, \ldots, a_1^{n_1}$, and note that this is an $m$-stackwise word that only differs from $a_i^{n_i}, \ldots, a_1^{n_1}$ by the additional $q$ added to the multiplicity of $a_j$. Set $k' = n_j + \cdots + n_1 + 1$. In the notation of Definition 4.7 we have $[b^p, a_{i-1}^{n_{i-1}}, \ldots, c^{q+r}, a_{j-1}^{n_{j-1}}, \ldots, a_1^{n_1}]_{k'} = x_j < x_{j+1} = [c^q, a_i^{n_i}, \ldots, a_1^{n_1}]_{k'}$ while the first $k' - 1$ symbols in each expression are identical. So in this case the output of the operator with $g = 0$ comes from an $m$-stackwise word that is less than $c^q, a_i^{n_i}, \ldots, a_1^{n_1}$ in the total word order, fitting scenario (e) for this operator.

Continue to suppose that $g = 0$ and again move $c^{q+r}$ to the right so that the initial segment of the word is $c^{q+r}, a_{j-1}^{n_{j-1}}, \ldots, a_1^{n_1}$, but now suppose this initial segment is not an $m$-stackwise word. If $j = 1$, then $y$ in minimal in $P$ and so we must have $q + r > m$ for $c^{q+r}$ to not be an $m$-stackwise word. But then the action vanishes as above, so assume $j \geq 2$. The length of $c^{q+r}, a_{j-1}^{n_{j-1}}, \ldots, a_1^{n_1}$ is less than the length of $c^q, a_i^{n_i}, \ldots, a_1^{n_1}$ since the former is at least missing $b^p$, and it has augmented stackwise form, so apply the SSTIH to $\langle c^{q+r}, a_{j-1}^{n_{j-1}}, \ldots, a_1^{n_1}\rangle.\{\emptyset, \ldots, \emptyset\}$. Since its word is not $m$-stackwise, either (a) or (c) or (e) from Theorem 6.1 applies. Both (a) and (c) result in a vanishing action, so for the remainder of this subcase we will assume that (e) applies.

Let $\mathbf{d}$ be an $m$-stackwise word from the expansion of $\langle c^{q+r}, a_{j-1}^{n_{j-1}}, \ldots, a_1^{n_1}\rangle.\{\emptyset, \ldots, \emptyset\}$ from (e) and note that $\mathbf{d} < c^{q+r}, a_{j-1}^{n_{j-1}}, \ldots, a_1^{n_1}$ in the total word order. Set $k' = n_j + \cdots + n_1$. Suppose that $R_{k'}(\mathbf{d}) \neq R_{k'}(c^{q+r}, a_{j-1}^{n_{j-1}}, \ldots, a_1^{n_1})$. Since $\mathbf{d} < c^{q+r}, a_{j-1}^{n_{j-1}}, \ldots, a_1^{n_1}$, we must have $R_{k'}(\mathbf{d}) < R_{k'}(c^{q+r}, a_{j-1}^{n_{j-1}}, \ldots, a_1^{n_1})$. Also, we have $R_{k'}(c^{q+r}, a_{j-1}^{n_{j-1}}, \ldots, a_1^{n_1}) = R_{k'}(c^q, a_i^{n_i}, \ldots, a_1^{n_1})$. Hence $\mathbf{d}$ is an $m$-stackwise word with $R_{k'}(\mathbf{d}) < R_{k'}(c^q, a_i^{n_i}, \ldots, a_1^{n_1})$, so using $h = k'$ we apply Lemma 6.2 to the expansion of $\langle a_{i-1}^{n_{i-1}}, \ldots, b^p, \mathbf{d}\rangle.\{\emptyset, \ldots, \emptyset\}$ to show (e) is satisfied for this kind of subword $\mathbf{d}$. Now suppose $R_{k'}(\mathbf{d}) = R_{k'}(c^{q+r}, a_{j-1}^{n_{j-1}}, \ldots, a_1^{n_1})$. The total color multiplicities of both words must be the same by Remark 3.2, so both words must end with $c^q$. This results in $\mathbf{d}$ and $c^{q+r}, a_{j-1}^{n_{j-1}}, \ldots, a_1^{n_1}$ being the same word, contradicting that $\mathbf{d} < c^{q+r}, a_{j-1}^{n_{j-1}}, \ldots, a_1^{n_1}$. Hence this case cannot happen. In the possible cases for $g = 0$, we have shown that (e) applies.

Now suppose that $g \geq 1$ so that the initial segment is $b^g, \mathbf{f}$. First suppose that $j = 1$ so that $y$ is a minimal element of $P$. Here $\mathbf{f} = f^s$, (i),(ii) and these colors are all distant to $c$. So the color group $b^g$ is the first appearance of either $b$ or $c$ in $b^g, \mathbf{f}$. No ideal in any $m$-multiset of ideals in the expansion to this point contains the minimal element $y$ and so no ideal can accept $x$ since it is the minimal element of color $b$ in $P$. Hence the result of the action vanishes.

Continue to assume $g \geq 1$ and now suppose $j \geq 2$. The element $y$ is minimal in the filter $P - K$ by our work in the paragraph introducing $K$, and $y$ has color $c$. Since $b$ and $c$ are adjacent, by AC all elements of these colors must be comparable. Hence $P - K$ cannot have a minimal element of color $b$, so $b^g, \mathbf{f}$ is not an $m$-stackwise word. This word is less than $c^q, a_i^{n_i}, \ldots, a_1^{n_1}$ in the total word order by length and has augmented stackwise form, so apply the SSTIH to write $\langle b^g, \mathbf{f} \rangle.\{\emptyset, \ldots, \emptyset\}$ as an integral sum of $m$-stackwise vectors whose words are less than or equal to $b^g, \mathbf{f}$ in the total word order. As before, only (a) or (c) or (e) can apply since $b^g, \mathbf{f}$ is not $m$-stackwise, and we will assume (e) since (a) and (c) result in vanishing actions.

Let $\mathbf{d}$ be an $m$-stackwise word resulting from this expansion. Set $k' = n_{j-1} + \cdots + n_1$. First suppose that $R_{k'}(\mathbf{d}) \neq R_{k'}(b^g, \mathbf{f})$. Then $R_{k'}(\mathbf{d}) < R_{k'}(b^g, \mathbf{f})$ since $\mathbf{d} < b^g, \mathbf{f}$. Also, we have $R_{k'}(b^g, \mathbf{f}) = R_{k'}(c^q, a_i^{n_i}, \ldots, a_1^{n_1})$. Hence $\mathbf{d}$ is an $m$-stackwise word with $R_{k'}(\mathbf{d}) < R_{k'}(c^q, a_i^{n_i}, \ldots, a_1^{n_1})$, so taking $h = k'$ we apply Lemma 6.2 to the expansion of $\langle a_{i-1}^{n_{i-1}}, \ldots, b^{p-g}, c^{q+r}, \mathbf{d} \rangle.\{\emptyset, \ldots, \emptyset\}$ to show (e) is satisfied for this kind of subword $\mathbf{d}$.

Finally suppose that $R_{k'}(\mathbf{d}) = R_{k'}(b^g, \mathbf{f}) = a_{j-1}^{n_{j-1}}, \ldots, a_1^{n_1}$. There are $g$ copies of $b$ from $b^g$ that must appear somewhere in $\mathbf{d}$ by Remark 3.2. We claim they must appear at $a_{j-1}$. They cannot appear in the initial segment $R_{k'}(\mathbf{d})$ since $R_{k'}(\mathbf{d}) = R_{k'}(b^g, \mathbf{f})$. The remaining colors appearing in $b^g, \mathbf{f}$ are adjacent (in the case of $b^g$) or distant (in the case of $f^s$, (i),(ii)) to $c$ and $y$ has color $c$. Since $\mathbf{d}$ and $b^g, \mathbf{f}$ have the same color census, the element $y$ is not in the ideal $L$ consisting of elements corresponding to $\mathbf{d}$. Since $b$ and $c$ are adjacent, we note by AC that $y$ would have to be in $L$ if any element of color $b$ greater than $y$ in the extension is in $L$. Hence these $g$ copies of $b$ cannot appear to the left of $a_{j-1}$ in $\mathbf{d}$. Thus they must appear at $a_{j-1}$, as claimed, and add to the total multiplicity at $a_{j-1}$. Set $k'' = k' + 1$ and note that $R_{k''}(\mathbf{d}) = a_{j-1}^{n_{j-1}+1}, a_{j-2}^{n_{j-2}}, \ldots, a_1^{n_1}$ and $R_{k''}(c^q, a_i^{n_i}, \ldots, a_1^{n_1}) = a_j, a_{j-1}^{n_{j-1}}, a_{j-2}^{n_{j-2}}, \ldots, a_1^{n_1}$. Both of these right segments grow well for $m$ by Lemma 4.3(b) since they are the initial segments of the $m$-stackwise words $\mathbf{d}$ and $a_i^{n_i}, \ldots, a_1^{n_1}$, respectively. Hence in the total word order of Definition 4.7 we have $R_{k''}(\mathbf{d}) < R_{k''}(c^q, a_i^{n_i}, \ldots, a_1^{n_1})$ since $[R_{k''}(\mathbf{d})]_\ell = [R_{k''}(c^q, a_i^{n_i}, \ldots, a_1^{n_1})]_\ell$ for all $1 \leq \ell \leq k'$ but $[R_{k''}(\mathbf{d})]_{k''} = x_{j-1} < x_j = [R_{k''}(c^q, a_i^{n_i}, \ldots, a_1^{n_1})]_{k''}$. Using $h = k''$ we apply Lemma 6.2 to the expansion of $\langle a_{i-1}^{n_{i-1}}, \ldots, b^{p-g}, c^{q+r}, \mathbf{d} \rangle.\{\emptyset, \ldots, \emptyset\}$ to show (e) is satisfied for this kind of subword $\mathbf{d}$. This completes the case where there is some element covered by $x$ which follows $y$ in the extension.

If there is no such element so that $y$ is the greatest element in the extension covered by $x$, then a similar proof proceeds as above without the considerations given to $f^s$ and the color groups corresponding to elements in groups (i) and (ii) in $P$. $\square$

Finishing the proofs of Lemmas 6.3, 6.4, and 6.5 completes the proof of Theorem 6.1.

## 8 Main result

Theorem 8.1 is the culmination of the work of Sections 4–7. Our goal has been to show that Γ-colored *d*-complete posets are standard, which we defined in Definition 3.3. Since Γ-colored minuscule posets are also Γ-colored *d*-complete, the results of these previous sections also apply to Γ-colored minuscule posets. We present representation theoretic applications in Section 9.

As noted in the Introduction, this work, including the following theorem, was jointly obtained with Robert A. Proctor.

**Theorem 8.1.** *Let P be a poset colored by a simply laced Dynkin diagram Γ. If P is Γ-colored d-complete or Γ-colored minuscule, then P is standard.*

*Proof.* Suppose that $P$ is Γ-colored *d*-complete. Let $m \geq 1$ and fix a linear extension of $P$. Proposition 4.6 shows that the *m*-stackwise vectors are linearly independent. Let $\langle b_k^{n_k}, \ldots, b_1^{n_1} \rangle.\{\emptyset, \ldots, \emptyset\}$ be an *m*-vector in $\mathcal{V}^m$. Theorem 6.1 shows that $\langle b_1^{n_1} \rangle.\{\emptyset, \ldots, \emptyset\}$ either vanishes or can be expanded into an integer linear sum of *m*-stackwise vectors, that each of these can be acted on with $\langle b_2^{n_2} \rangle$ to either vanish or obtain a new integer linear sum of *m*-stackwise vectors, and so on. This expansion results in the original *m*-vector being zero or written as an integer linear sum of *m*-stackwise vectors. Hence $P$ is standard. The result also holds for Γ-colored minuscule posets since they are Γ-colored *d*-complete. □

**Remark 8.2.** The definitions of Γ-colored *d*-complete and Γ-colored minuscule posets need not be restricted to simply laced Γ. The classifications of each kind given in [Str2] and [Str3] were for general Dynkin diagrams. We conjecture that Theorem 8.1 also holds for Dynkin diagrams that are not simply laced. The main hurdle in extending this result is developing new operator identities from Section 5 that hold when intervals in *P* may consist of a chain of three elements whose bottom and top elements have the same color. Such intervals exist in the general case but do not exist in the simply laced case.

## 9 Implications for minuscule representations and the Seshadri basis

In this section we use Lie theoretic terminology such as Kac–Moody algebra, representation, module, weight, weight diagram, Weyl group, etc; consult [Hum] and [Kac] and [Kum] for definitions as needed.

We now fix an order $a_1, \ldots, a_n$ of Γ. This makes the integers $\{\theta_{ij}\}_{a_i,a_j \in \Gamma}$ from Remark 2.1(b) the entries of a generalized Cartan matrix (GCM). The Lie algebras associated to Γ and its GCM that we are concerned with are the Kac–Moody algebra $\mathfrak{g}$ and several of its subalgebras, including its derived subalgebra $\mathfrak{g}' = [\mathfrak{g}, \mathfrak{g}]$ (where $[\cdot, \cdot]$ is the Lie bracket for $\mathfrak{g}$), its Cartan subalgebra $\mathfrak{h}$, its positive nilpotent subalgebra $\mathfrak{n}$, and its positive Borel subalgebra $\mathfrak{b} = \mathfrak{h} + \mathfrak{n}$. We define the *Cartan derived subalgebra* and the *Borel derived subalgebra* as $\mathfrak{h}' = \mathfrak{h} \cap \mathfrak{g}'$ and $\mathfrak{b}' = \mathfrak{b} \cap \mathfrak{g}'$, respectively. The algebras $\mathfrak{h}'$, $\mathfrak{n}$, $\mathfrak{b}'$, and $\mathfrak{g}'$ can be realized via generators taken from $\{x_i, y_i, h_i\}_{a_i \in \Gamma}$ and relations taken from the following list.

(XX) For all $a_i, a_j \in \Gamma$ such that $i \neq j$ we have $\underbrace{[x_i, [\ldots, [x_i, x_j] \ldots]]}_{1-\theta_{ji} \text{ times}} = 0$, (YY) For all $a_i, a_j \in \Gamma$ such that $i \neq j$ we have $[\underbrace{y_i, [\ldots, [y_i, y_j]\ldots]}_{1-\theta_{ji} \text{ times}}] = 0$,

(HH) For all $a_i, a_j \in \Gamma$ we have $[h_j, h_i] = 0$,

(HX) For all $a_i, a_j \in \Gamma$ we have $[h_j, x_i] = \theta_{ij} x_i$,

(HY) For all $a_i, a_j \in \Gamma$ we have $[h_j, y_i] = -\theta_{ij} y_i$,

(XY) For all $a_i, a_j \in \Gamma$ we have $[x_i, y_j] = \delta_{ij} h_i$, where $\delta_{ij}$ is the Kronecker delta.

The algebra $\mathfrak{h}'$ is generated by $\{h_i\}_{a_i \in \Gamma}$ subject to the relation HH. The algebra $\mathfrak{n}$ is generated by $\{x_i\}_{a_i \in \Gamma}$ subject to the relation XX. The algebra $\mathfrak{b}'$ is generated by $\{x_i, h_i\}_{a_i \in \Gamma}$ subject to the relations XX, HH, and HX. We note that $\mathfrak{b}' = \mathfrak{h}' + \mathfrak{n}$. The algebra $\mathfrak{g}'$ is generated by $\{x_i, y_i, h_i\}_{a_i \in \Gamma}$ subject to all of the relations above. If the GCM for $\Gamma$ is invertible, then $\mathfrak{h}$, $\mathfrak{b}$, and $\mathfrak{g}$ are equal to their derived counterparts.

We use the color raising and lowering operators $\{X_i\}_{a_i \in \Gamma}$ and $\{Y_i\}_{a_i \in \Gamma}$ defined on $V = \langle \mathcal{FI}(P) \rangle$ in Section 2 to construct representations of the Lie algebras $\mathfrak{b}'$ and $\mathfrak{g}'$.

**Definition 9.1.** Let $P$ be a $\Gamma$-colored poset and let $\mathfrak{L}$ be the Lie algebra $\mathfrak{b}'$ or $\mathfrak{g}'$. We say that $V$ *carries a representation of* $\mathfrak{L}$ if there are operators $\{H_i\}_{a_i \in \Gamma}$ diagonal on the basis $\{\langle F, I \rangle\}_{(F,I) \in \mathcal{FI}(P)}$ of $V$ such that the operators $\{X_i, H_i\}_{a_i \in \Gamma}$ (if $\mathfrak{L} = \mathfrak{b}'$) or $\{X_i, Y_i, H_i\}_{a_i \in \Gamma}$ (if $\mathfrak{L} = \mathfrak{g}'$) satisfy the defining relations for $\mathfrak{L}$ under the commutator bracket $[A, B] = AB - BA$ on $\text{End}(V)$.

If $V$ carries a representation of $\mathfrak{b}'$ or $\mathfrak{g}'$, then the map $x_i \mapsto X_i$, $y_i \mapsto Y_i$ (if applicable), $h_i \mapsto H_i$ induces a Lie algebra homomorphism from $\mathfrak{b}'$ or $\mathfrak{g}'$ to the subspace of $\text{End}(V)$ generated by $\{X_i, H_i\}_{a_i \in \Gamma}$ or $\{X_i, Y_i, H_i\}_{a_i \in \Gamma}$ under the commutator bracket. Hence $V$ is a representation of $\mathfrak{b}'$ or $\mathfrak{g}'$. In general many choices for the operators $\{H_i\}_{a_i \in \Gamma}$ will work; see [Str1, §6].

As noted in the Introduction, the most notable representations in [Str1] share characteristics with the classic minuscule representations of simple Lie algebras. We continue to emphasize ideals and denote the action $X_i.\langle F, I \rangle$ as $X_i.I$, and similarly for the other operators.

**Definition 9.2.** Let $P$ be a $\Gamma$-colored poset.
  (a) A representation of $\mathfrak{b}'$ carried by $V$ is *upper P-minuscule* if
      (i) We have $X_i^2 = 0$ for all $a_i \in \Gamma$,
      (ii) The eigenvalues for $H_i$ for all $a_i \in \Gamma$ are contained in $\{-1, 0, 1, 2, 3, \ldots\}$, and
      (iii) For every split basis vector $\langle F, I \rangle$ and every $a_i \in \Gamma$ we have $H_i.I = -I$ if and only if $F = P - I$ contains a minimal element of color $a_i$, i.e. if and only if $X_i.I \neq 0$.
  (b) A representation of $\mathfrak{g}'$ carried by $V$ is *P-minuscule* if the eigenvalues for $H_i$ for all $a_i \in \Gamma$ are contained in $\{-1, 0, 1\}$.

We state the main result of [Str1] which characterizes these representations. We continue to work only with finite posets in this paper, though the result restated here includes infinite locally finite posets as well.

**Theorem 9.3** (Thm. 38 of [Str1]). *Let $P$ be a $\Gamma$-colored poset.*
  (a) *The vector space $V$ carries an upper P-minuscule representation of $\mathfrak{b}'$ if and only if $P$ is a $\Gamma$-colored d-complete poset.*

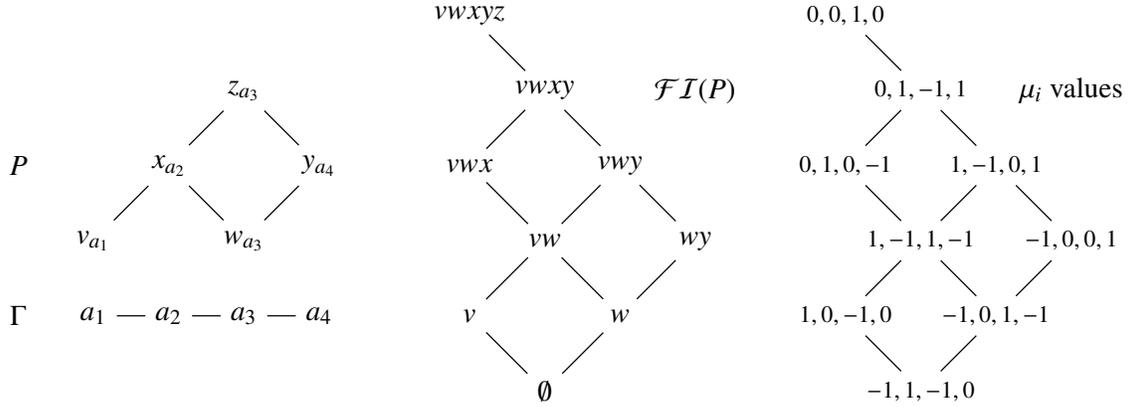

Figure 9.1: A relabeling and continuation of Figure 3.1. On the right we have added values for the $\mu$-diagonal operators on $\mathcal{FI}(P)$. For each split $(F, I) \in \mathcal{FI}(P)$ its $\mu_i$ values are displayed in the corresponding place on the right in the format $\mu_1(I), \mu_2(I), \mu_3(I), \mu_4(I)$.

*(b) The vector space V carries a P-minuscule representation of $\mathfrak{g}'$ if and only if P is a $\Gamma$-colored minuscule poset.*

The operators $\{H_i\}_{a_i \in \Gamma}$ are unique for upper $P$-minuscule and $P$-minuscule representations; see [Str1, Prop. 33, Thm. 35]. We describe these operators in that setting. Let $P$ be $\Gamma$-colored $d$-complete with corresponding upper $P$-minuscule representation $V$. Then $P$ has a unique maximal element (see [Str3, Lem. 3.3]). Let $a_j$ be the color of this unique maximal element. For the maximal ideal $P$ of $\mathcal{FI}(P)$ we define $\mu_j(P) = 1$ and $\mu_i(P) = 0$ whenever $i \neq j$. The results of [Str1, §6] show that this choice uniquely determines a collection of functions $\{\mu_i\}_{a_i \in \Gamma}$ that produce the operators $H_i.I = \mu_i(I)I$ we desire. Adapting Equation (3) of [Str1, pg. 15] to our situation, we obtain the following recursive description of the functions $\{\mu_i\}_{a_i \in \Gamma}$. Suppose $I \to I'$ in $\mathcal{FI}(P)$ and $I' - I$ consists of an element of color $a_k$. If $\mu_i(I')$ is known, then $\mu_i(I) = \mu_i(I') - 2\delta_{ik} + \sum_{a_\ell \sim a_k} \delta_{i\ell}$ where $\delta$ is the Kronecker delta. It can be seen that this produces the $\mu$-diagonal operators of [Str1, §7] by starting with $P$ and working our way down through $\mathcal{FI}(P)$.

**Example 9.4.** Consider again the $\Gamma$-colored $d$-complete poset $P$ and distributive lattice $\mathcal{FI}(P)$ from Example 3.1. Here we work out the values taken by $\{\mu_i\}_{1 \leq i \leq 4}$ and display our results in Figure 9.1. We note that computing the $\mu$-diagonal operators in this way amounts to playing the Numbers Game of Mozes [Moz] starting from the way the operators $\{\mu_i\}$ were defined on the maximal ideal $P$.

Let $P$ be a $\Gamma$-colored $d$-complete poset with upper $P$-minuscule representation $V$. Consider the $m^{\text{th}}$ tensor power $V^{\otimes m}$ of $V$. Its basis $\{I_1 \otimes \cdots \otimes I_m \mid (F_j, I_j) \in \mathcal{FI}(P)\}$ can be realized as the ordered $m$-tuples of ideals $(I_1, \ldots, I_m)$ of Section 5. Then $V^{\otimes m}$ is a representation of $\mathfrak{b}'$ under the action $b.(I_1, \ldots, I_m) = \sum_{k=1}^{m}(I_1, \ldots, b.I_k, \ldots, I_m)$ for $b \in \mathfrak{b}'$. When $b = x_i$ for some $a_i \in \Gamma$ we recover how the operator $X_i$ was defined in Section 5. Now consider the $m^{\text{th}}$ symmetric power $S^m V$ of $V$ which consists of the degree $m$ elements of the tensor algebra of $V$ after quotienting by the two sided ideal generated by elements of the form $I_r \otimes I_s - I_s \otimes I_r$. For our purposes this allows positions within each $m$-tuple of ideals to commute, producing unordered $m$-multisets of ideals $\{I_1, \ldots, I_m\}$. The action is defined in the same way, recovering how the operator $X_i$ was defined in Section 3. These comments hold for $\Gamma$-colored minuscule posets $P$ and $P$-minuscule representations $V$ of $\mathfrak{g}'$ also.

Theorem 8.1 showed $P$ is standard, so the $m$-stackwise vectors are a basis for $\text{span}_{\mathbb{C}}(\mathcal{V}^m) \subseteq S^m V$. Proctor plans to show in a future work (personal communication) that when $m = 1$, the $\mathfrak{b}'$-module $V$ is isomorphic to the Demazure module for an integrable highest weight module $V(\lambda)$. This module is produced using the Weyl group element $w$ corresponding to the playing of the Numbers Game that produced the weight diagram of $V$. This weight diagram is a "minuscule upper portion" of the weight diagram of $V(\lambda)$. For instance, we are working in type $A_4$ in Example 9.4 with $\lambda$ being the fundamental weight $\omega_3$ and $w = s_1 s_3 s_2 s_4 s_3$, where $\{s_1, s_2, s_3, s_4\}$ are the Weyl group generators. So $\text{span}_{\mathbb{C}}(\mathcal{V}^m)$ is a $\mathfrak{b}'$-module sitting inside of the $m^{\text{th}}$ symmetric power $S^m V$ of this Demazure module $V$.

We turn our attention to semisimple Lie algebras and $\Gamma$-colored minuscule posets. Note that $\Gamma$ has finite Lie type and so $\mathfrak{g}' = \mathfrak{g}$ and $\mathfrak{h}' = \mathfrak{h}$ in this setting. Let $\{\omega_1, \ldots, \omega_n\}$ be the fundamental weights associated to $\Gamma$. We have $\omega_i(h_j) = \delta_{ij}$ for $1 \leq i, j \leq n$. The *weight* of the split vector $\langle F, I \rangle$ in $V$ is

$$\text{wt}(I) = \sum_{i=1}^{n} \mu_i(I) \omega_i.$$

Thus $\text{wt}(P) = \omega_j$, where $a_j$ is the color of the maximal element of $P$. We denote this *highest weight* by $\lambda$. As noted in Example 9.4, computing the rest of the weights amounts to playing the Numbers Game of Mozes which produces the weight diagram of the representation. The recursive formula given prior to Example 9.4 becomes simply $\text{wt}(I) = \text{wt}(I') - \alpha_k$ when converted to weights and simple roots. That is, removing a single element of color $a_k$ from an ideal $I'$ corresponds to subtracting the simple root $\alpha_k$ from the weight of the corresponding vector, so the weight lattice of $V(\lambda)$ and the split lattice $\mathcal{FI}(P)$ have the same poset structure.

When $b = h_i$ for some $a_i \in \Gamma$ we get $h_i.(I_1, \ldots, I_m) = (\mu_i(I_1) + \cdots + \mu_i(I_m))(I_1, \ldots, I_m)$. Hence this action produces weight for $(I_1, \ldots, I_m)$ in $V^{\otimes m}$ of

$$\text{wt}(I_1, \ldots, I_m) = \sum_{k=1}^{m} \text{wt}(I_k) = \sum_{k=1}^{m} \sum_{i=1}^{n} \mu_i(I_k) \omega_i.$$

This computation also holds for $\{I_1, \ldots, I_m\} \in S^m V$. The weight of the highest weight vector $(P, \ldots, P)$ in $V^{\otimes m}$ and $\{P, \ldots, P\}$ in $S^m V$ is thus $m\omega_j = m\lambda$.

The $\Gamma$-colored minuscule posets were classified in [Str3]; Figure 8.1 of that paper displays all finite $\Gamma$-colored posets. This classification was produced from the $\Gamma$-colored minuscule defining axioms, but the list of such finite posets is exactly the list of colored minuscule posets of [Pr1, Thm. 11]. These latter posets were obtained from the weight diagrams of the classic minuscule representations of simple Lie algebras. Such representations are indexed by their highest weights; a list of all minuscule weights is given in [Bou, Ch. VIII, §7.3]. There is one minuscule representation and one $\Gamma$-colored minuscule poset for each minuscule weight. Each colored minuscule poset in [Pr1] was obtained as the poset of join irreducible elements of the corresponding minuscule weight diagram, which itself was produced by playing the Numbers Game. The coloring of those posets in [Pr1, Thm. 11] was given after identifying elements of the colored minuscule posets with certain coroots; see also [Str3, Thm. 9.2].

Let $P$ be a $\Gamma$-colored minuscule poset. The highest weight $\lambda = \omega_j$ is a minuscule weight. The corresponding $P$-minuscule representation of $\mathfrak{g}$ is precisely the irreducible representation $V(\lambda)$, that is, $V \simeq V(\lambda)$. Since $P$ is $\Gamma$-colored minuscule, we note that the order dual $P^*$ is also $\Gamma$-colored $d$-complete. The $\{Y_i\}_{a_i \in \Gamma}$ operators are dualized versions of the $\{X_i\}_{a_i \in \Gamma}$ operators and hence all of the results of this paper for $\Gamma$-colored $d$-complete posets and $\{X_i\}_{a_i \in \Gamma}$ operators working upward from $\{\emptyset, \ldots, \emptyset\}$ hold for $\{Y_i\}_{a_i \in \Gamma}$ operators working downward from $\{P, \ldots, P\}$. We used $m$-flags of ideals in the upward case, while $m$-flags of filters are used in the downward case. To distinguish the two scenarios, we refer to *upward* structures in the former case and *downward* structures in the latter case and use subscripts $U$ and $D$ similarly; for example, Theorem 8.1 shows that the upward $m$-stackwise vectors form a basis for the $\mathfrak{g}$-module $\mathrm{span}_{\mathbb{C}}(\mathcal{V}^m)_U$ over the integers. By dualizing our work, this theorem also shows that the downward $m$-stackwise vectors form a basis for $\mathrm{span}_{\mathbb{C}}(\mathcal{V}^m)_D$ over the integers. Both correspond to $m$-multichains in the weight lattice of $V(\lambda)$ since this weight lattice has the same poset structure as the lattice $\mathcal{FI}(P)$. As noted above, the highest weight vector $\{P, \ldots, P\}$ has weight $m\lambda$. The submodule of $S^m V$ generated by this downward action of the $\{Y_i\}_{a_i \in \Gamma}$ operators on $\{P, \ldots, P\}$ is isomorphic to $V(m\lambda)$, i.e. $\mathrm{span}_{\mathbb{C}}(\mathcal{V}^m)_D \simeq V(m\lambda)$.

We summarize these observations in the following theorem. As noted in the Introduction, this work, including this theorem, was jointly obtained with Robert A. Proctor.

**Theorem 9.5.** *Let $P$ be a $\Gamma$-colored minuscule poset for simply laced $\Gamma$ and let $m \geq 1$. Let $\lambda = \omega_j$, where $a_j$ is the color of the unique maximal element of $P$. Then $\lambda$ is a minuscule weight and $\mathrm{span}_{\mathbb{C}}(\mathcal{V}^m)_D$ is isomorphic to the irreducible representation $V(m\lambda)$ of highest weight $m\lambda$. The downward $m$-stackwise vectors form a basis for $\mathrm{span}_{\mathbb{C}}(\mathcal{V}^m)_D$ and are in natural one-to-one correspondence with $m$-multichains in the weight lattice of $V(\lambda)$. All $m$-vectors in $\mathcal{V}_D^m$ can be expanded as integer linear sums of the basis of downward $m$-stackwise vectors.*

This indexing of downward $m$-stackwise basis vectors with such $m$-multichains emulates Seshadri's indexing in the simply laced cases of Theorem 1.1 above.

## Acknowledgements

I would like to thank Robert A. Proctor for many contributions to this paper. First, for his original 1993 manuscript, without which this paper would not be possible. And second, for thoroughly reading through most of an earlier version of this manuscript and offering many helpful corrections and comments, including recommending important corrections in Sections 4, 5, and 6. I would also like to thank Sam Jeralds for providing helpful feedback on Sections 1 and 9.